\documentclass[reqno,12pt]{amsart}

\usepackage{    amsmath,
                amsfonts,
                amssymb,
                amsthm,
                amsxtra,
                verbatim,
}

\usepackage[dvips]{graphics}

\InputIfFileExists{diagrams.tex}{}{%
\IfFileExists{diagrams.sty}{\usepackage{diagrams}}
{\message{Without diagrams you can not process this file!}}}
\diagramstyle[height=2em,balance,righteqno,PostScript=dvips,%
dpi=600]        %%     to print with 600 dots per inch

\def\rhaha{\raise.24ex\hbox{$\rightharpoonup$}\kern-1em\lower.24ex\hbox{$\rightharpoondown$}}%
\def\lhaha{\raise.24ex\hbox{$\leftharpoonup$}\kern-1em\lower.24ex\hbox{$\leftharpoondown$}}%
\def\dhaha{\downharpoonleft\kern-.22em\downharpoonright\kern.02em}%
\def\uhaha{\upharpoonleft\kern-.22em\upharpoonright\kern.02em}% -.21em
\newarrowhead{twoharpoons}\rhaha\lhaha\dhaha\uhaha

\newarrow{Id}===={twoharpoons}
\newarrow{Epi}----{triangle}
\newarrow{Twoar}===={=>}
\newarrow{Mono}{boldhook}---{->}
\newarrow{TTo}----{->}

\newcommand\CC{{\mathbb C}}
\newcommand\FF{{\mathbb F}}

\newcommand\RR{{\mathbb R}}
\newcommand\QQ{{\mathbb Q}}
\newcommand\ZZ{{\mathbb Z}}

\newcommand{\cc}{{\mathcal C}}

\newcommand{\cf}{{\mathcal F}}
\newcommand{\ch}{{\mathcal H}}
\newcommand{\cj}{{\mathcal J}}
\newcommand{\co}{{\mathcal O}}
\newcommand{\Top}{{\mathcal T}}
\newcommand{\cu}{{\mathcal U}}
\newcommand{\CU}{{\mathcal U}}
\newcommand{\cv}{{\mathcal V}}
\newcommand{\cw}{{\mathcal W}}
\newcommand{\cx}{{\mathcal X}}
\newcommand{\cy}{{\mathcal Y}}
\newcommand{\cz}{{\mathcal Z}}

\newcommand{\SSS}{{\mathfrak S}}

\newcommand{\Dbc}{D^{b,c}}

\newcommand{\0}{\phantom0}

\newcommand{\bull}{{\scriptscriptstyle\bullet}}

\newcommand{\nquad}{\!\!\!\!\!\!}

\newcommand{\n}[1]{\nobreakdash-\hspace{0pt}}
\newcommand{\kd}[1]{$\Bbbk$\nobreakdash-\hspace{0pt}}

\let\bs\backslash

\let\dlgn\boxtimes

\let\tens\otimes
\let\ttt\textstyle
\let\und\underline
\let\xra\xrightarrow

\DeclareMathOperator\assoc{assoc}
\DeclareMathOperator\Aut{Aut}
\DeclareMathOperator\Can{Can}
\DeclareMathOperator\coass{coass}
\DeclareMathOperator\coher{coher}
\DeclareMathOperator\Coher{Coher}
\DeclareMathOperator\const{const}

\DeclareMathOperator\Gr{Gr}
\DeclareMathOperator\Hom{Hom}
\DeclareMathOperator\id{id}
\DeclareMathOperator\Id{Id}

\DeclareMathOperator\Ker{Ker}

\DeclareMathOperator\Mor{Mor}
\DeclareMathOperator\Ob{Ob}
\DeclareMathOperator\pr{pr}

\DeclareMathOperator\Resto{Res}

\DeclareMathOperator\sets{-sets-}
\newcommand{\Sh}[1]{Sh(#1)}

\DeclareMathOperator\supp{supp}

\DeclareMathOperator\Vect{Vect}

\newcommand{\corref}[1]{Corollary~\ref{#1}}

\newcommand{\figref}[1]{Figure~\ref{#1}}
\newcommand{\lemref}[1]{Lemma~\ref{#1}}
\newcommand{\propref}[1]{Proposition~\ref{#1}}

\newcommand{\secref}[1]{Section~\ref{#1}}

\newcommand{\edc}[2]{\genfrac{}{}{}1{#2}{#1}}
\newcommand{\eDc}[2]{\genfrac{}{}{}{}{#2}{#1}}
\newcommand{\dec}[2]{\genfrac{}{}{}1{#1}{#2}}
\newcommand{\Dec}[2]{\genfrac{}{}{}{}{#1}{#2}}

\newtheoremstyle{boldhead}%     name
{\topsep}%                      abovespace
{\topsep}%                      belowspace
{\slshape}%         bodyfont
{}%             indentation=noindent
{\bfseries}%            headfont
{.}%                            headpunctuation
{ }%                headspace=interword space
{\thmname{#1}\thmnumber{ #2}\thmnote{ (#3)}}%   custom head specification

\newtheoremstyle{boldremark}%   name
{\topsep}%                      abovespace
{\topsep}%                      belowspace
{\upshape}%         bodyfont
{}%             indentation=noindent
{\bfseries}%            headfont
{.}%                            headpunctuation
{ }%                headspace=interword space
{\thmname{#1}\thmnumber{ #2}\thmnote{ (#3)}}%   custom head specification

\swapnumbers

\theoremstyle{boldhead}

\newtheorem{theorem}[subsubsection]{Theorem}
\newtheorem*{subtheorem}{Theorem}
\newtheorem{corollary}[subsubsection]{Corollary}
\newtheorem{lemma}[subsubsection]{Lemma}

\newtheorem{proposition}[subsubsection]{Proposition}

\theoremstyle{boldremark}

\newtheorem{definition}[subsubsection]{Definition}

\usepackage{t-angles} % for TPIC drivers such as dviwin, xdvi, dvips
\makeatletter
\def\@seccntformat#1{\csname the#1\endcsname.\quad}
\renewcommand\section{\@startsection {section}{1}{\z@}%
                                   {-3.5ex \@plus -1ex \@minus -.2ex}%
                                   {2.3ex \@plus.2ex}%
                                   {\normalfont\large\bfseries}}
\renewcommand\subsection{\@startsection{subsection}{2}{\z@}%
                                     {-3.25ex\@plus -1ex \@minus -.2ex}%
                                     {1.5ex \@plus .2ex}%
                                     {\normalfont\normalsize\bfseries}}
\renewcommand\subsubsection{\@startsection{subsubsection}{3}{\z@}%
                        {3.25ex plus 1ex minus .2ex}{-.5em}%
                        {\normalfont\normalsize\bfseries}}
\makeatother
\begin{document}
\title{The triangulated Hopf category $n_+SL(2)$}
\author{Volodymyr Lyubashenko}
\address{
Institute of Mathematics\\
National Academy of Sciences of Ukraine\\
3, Teresh\-chen\-kiv\-ska st.\\
Kyiv-4, 01601\\
Ukraine
}
\email{lub@imath.kiev.ua}
\date{math.QA/9904108; Edited March 28, 2001; Printed \today}
\begin{abstract}
We discuss an example of a triangulated Hopf category related to SL(2).
It is an equivariant derived category equipped with multiplication
and comultiplication functors and structure isomorphisms.
We prove some coherence equations for structure isomorphisms.
In particular, the Hopf category is monoidal.
\end{abstract}
\thanks{Research was supported in part by National Science Foundation
grant 530666.}
%\thanks{2000 \textit{Mathematics Subject Classification.}
%Primary 16W30, 18D15, 17B37; Secondary 18D35.}
\maketitle
\newtheorem*{acknowledgements}{Acknowledgements}
\newtheorem*{notations}{Notations}
%\end{comment}

\newlength{\texthigh}
\setlength{\texthigh}{\textheight}
\addtolength{\texthigh}{-0.12pt}
\newlength{\textwids}
\setlength{\textwids}{\textwidth}
\addtolength{\textwids}{-0.12pt}

\providecommand{\zhe}{\mathsf{X\mkern-8.2mu I\mkern3mu}}
\allowdisplaybreaks[1]

\section{Introduction}
Crane and Frenkel proposed a notion of a Hopf category
\cite{CraneFrenk:4-dim}. It was motivated by Lusztig's approach to
quantum groups -- his theory of canonical bases. In particular, Lusztig
obtains braided deformations $U_q\mathfrak n_+$ of universal enveloping
algebras $U\mathfrak n_+$ for some nilpotent Lie algebras
$\mathfrak n_+$ together with canonical bases of these braided Hopf
algebras \cite{Lus:canonical,Lus:quiver,Lus:book}.  Elements of the
canonical basis are identified with isomorphism classes of simple
perverse sheaves -- certain objects of equivariant derived categories.
They are contained in the subcategory of semisimple complexes. One of
the proposals of Crane and Frenkel is to study this category rather
than its Grothendieck ring $U_q\mathfrak n_+$. Conjectural properties
of this category were collected into a system of axioms of a Hopf
category, equipped with functors of multiplication and
comultiplication, isomorphisms of associativity, coassociativity and
coherence which satisfy four equations \cite{CraneFrenk:4-dim}.  A
mathematical framework and some examples were provided by
Neuchl~\cite{Neuchl-rtHc}.

Crane and Frenkel \cite{CraneFrenk:4-dim} gave an example of a Hopf
category resembling the semisimple category encountered in Lusztig's
theory corresponding to one-dimen\-sio\-nal Lie algebra $\mathfrak n_+$
-- nilpotent subalgebra of $\mathfrak{sl}(2)$. We want to discuss an
example of a related notion -- triangulated Hopf category -- the whole
equivariant derived category equipped with multiplication and
comultiplication functors and structure isomorphisms.  In particular,
it is a monoidal category.  New feature of coherence is that additive
relations of \cite{CraneFrenk:4-dim} are replaced  with distinguished
triangles.  This new structure does not induce a Hopf category
structure of Crane and Frenkel on the subcategory of semi-simple
complexes. The missing component is a consistent choice of splitting of
splittable triangles.  Verification of some of the consistency
equations is still an open question.

To give more details let us first recall some braided Hopf algebra
$H$. As an algebra $H$ is the algebra of polynomials of one variable
over $R=\ZZ[q,q^{-1}]$. More precisely, $H\subset\QQ(q)[x]$ is an
$R$\n-submodule spanned by the elements
\[ y_n = \frac{x^n}{(n)_{q^{-2}}!}, \qquad n\ge0, \]
\[(n)_{q^{-2}}! = \prod_{m=1}^n (m)_{q^{-2}}, \qquad (m)_{q^{-2}} =
\frac{1-q^{-2m}}{1-q^{-2}}. \]
The basis $(y_n)_{n\ge0}$ is called a canonical basis.

Multiplication table in this basis is
\[ y_k\cdot y_l = \binom{k+l}k_{q^{-2}} y_{k+l},\qquad
\binom{k+l}k_{q^{-2}} = \frac{(k+l)_{q^{-2}}!}{(k)_{q^{-2}}!
(l)_{q^{-2}}!} \in R. \]
Comultiplication by definition is
\[ y_n= \sum_{k+l=n} y_k\tens y_l. \]
These operations make $H$ into a $\ZZ$\n-graded $R$\n-algebra and
coalgebra. We equip the category $\cc$ of $\ZZ$\n-graded
$R$\n-modules with the braiding
\[ c: M\tens_R N \to N\tens_R M, \qquad c = \sum_{k,l} c_{k,l}, \]
\[
\begin{tanglec}
\object{k}\Step\object{l} \\
\xx \\
\object{l}\Step\object{k}
\end{tanglec}
= c_{k,l}: M_k\tens N_l \to N_l\tens M_k, \qquad c_{k,l}=q^{-2kl}, \]
where $M=\sum_{k\in\ZZ}M_k$, $N=\sum_{l\in\ZZ}N_l$ are graded
$R$\n-modules. Then $H$ is a Hopf algebra in a braided category $\cc$ as
defined in \cite{Ma:bg}. Using graphical notations
\[ m=\sum_{k,l}m_{k,l}, \qquad
\vstretch 60
\begin{tanglec}
\object{k}\Step\object{l} \\
\d\dd \\
\s \\
\object{k+l}
\end{tanglec}
=m_{k,l}: H_k\tens_R H_l \to H_{k+l}, \]
\[ \Delta=\sum_{k,l}\Delta_{k,l}, \qquad
\vstretch 60
\begin{tanglec}
\object{k+l} \\
\n \\
\dd\d \\
\object{k}\Step\object{l}
\end{tanglec}
=\Delta_{k,l}: H_{k+l} \to H_k\tens_R H_l, \]
we can write the coherence axiom as an equation
\begin{equation}\label{HnHmHpHq}
\vstretch 120
\begin{tanglec}
\nodeld{n}\Step\noderd{m}\\
\nw1\nodeld{n+m}\node\ne1\\
\n\\
\nodelu{p}\sw1\se1\noderu{q}
\end{tanglec}
\quad\; = \quad
\sum_{\substack{i+j=n \\ k+l=m \\ i+k=p \\ j+l=q}} \quad
\vstretch 72
\begin{tanglec}
\nodeld{n}\step[3]\noderd{m}\\
\nodeld{i}\n\noderd{j}\step[3]\nodeld{k}\n\noderd{l}\\
\sw1\se1\step\sw1\se1\\
\id\Step\hxx\Step\id\\
\nw1\ne1\step\nw1\ne1\\
\nodelu{p}\s\step[3]\s\noderu{q}
\end{tanglec}
\quad: H_n\tens_R H_m \to H_p\tens_R H_q
\end{equation}
for all $n,m,p,q\in\ZZ_{\ge0}$ such that $n+m=p+q$.

This algebra was obtained by Lusztig\cite{Lus:book} from the
following setup:

$\ch_{n_1,\dots,n_k}$ are $\CC$\n-linear categories, depending
symmetrically on parameters $n_1$,\dots,$n_k\in\ZZ_{\ge0}$,
$k\in\ZZ_{\ge0}$;

$m_{k,l}:\ch_{k,l}\to\ch_{k+l}$, $\Delta_{k,l}:\ch_{k+l}\to\ch_{k,l}$
are $\CC$\n-linear functors of
multiplication and comultiplication;

$c_{k,l}:\ch_{k,l}\to\ch_{l,k}$ are braiding functors;

there are associativity isomorphisms
\begin{diagram}
\ch_{k,l,n} & \rTTo^{1\dlgn m_{l,n}} & \ch_{k,l+n} \\
\dTTo<{m_{k,l}\dlgn1} & \simeq & \dTTo>{m_{k,l+n}} \\
\ch_{k+l,n} & \rTTo^{m_{k+l,n}} & \ch_{k+l+n}
\end{diagram}
where the meaning of $\dlgn$ will be specified further;

there are similar coassociativity isomorphisms.

The category $\ch_{n_1,\dots,n_k}$ is
$D^{b,c}_{GL(n_1)\times\dots\times GL(n_k)}(pt)$ -- the bounded
constructible equivariant derived category of a point. It has a
distinguished object $Y_{n_1,\dots,n_k}$ -- the constant sheaf,
which is the complex
\[ \dots \to 0 \to 0 \to \CC \to 0 \to 0 \to \dots \]
concentrated in degree 0. It turns out that the collection
$(Y_{n_1,\dots,n_k})$ is closed under multiplication and
comultiplication (up to coefficients which are graded vector spaces):
\[ m_{k,l}(Y_{k,l}) \simeq H^\bull(\Gr_k^{k+l}(\CC),\CC)\tens_\CC
Y_{k+l}. \]
The coefficient vector space here is de Rham
cohomology of the Grassmannian $\Gr_k^{k+l}(\CC)$ -- manifold of
$k$\n-dimensional subspaces of a $k+l$\n-dimensional space.
Cohomology is concentrated in even degrees and the Betti numbers
\[ \beta_i = \dim_\CC H^{2i}(\Gr_k^{k+l}(\CC),\CC) \]
are coefficients of the expansion of a $q$\n-binomial coefficient
in powers of $q^{-2}$:
\[ \binom{k+l}k_{q^{-2}} = \sum_{i\ge0} \beta_i q^{-2i}. \]
Replacing the degree with the power of $q$ we get the multiplication
table for the canonical basis $(y_k)$. Comultiplication law is
recovered from
\[ \Delta_{k,l}(Y_{k+l}) \simeq Y_{k,l}. \]

The braiding functor
\[ c_{k,l} = \bigl( D_{GL(k)\times GL(l)}(pt) \simeq D_{GL(l)\times
GL(k)}(pt) \rTTo^{[-2kl]} D_{GL(l)\times GL(k)}(pt) \bigr) \]
is essentially the degree shift by $-2kl$. It translates into
multiplication by $q^{-2kl}$ for the braiding in algebra setting.

In the present paper we shall discuss coherence at the category
level. If one replaces linear mappings in equation~\eqref{HnHmHpHq}
with functors and $\sum$ with $\oplus$
the equation fails: the left and the right hand side functors
$\ch_{n,m}\to\ch_{p,q}$ are, in general, not isomorphic. (Restricted
to $Y_{n,m}$ they give, however, isomorphic results.) One of the results
of the present paper is the following. Value of the left hand
side functor on an object $X$ of $\ch_{n,m}$ is a repeated extension
of values on $X$ of summands in the right hand side in the sense
of distinguished triangles. Precise analogy is as follows: a sheaf
$S$ on a topological space $W$ is an
extension of its quotient-sheaf $S_F$ supported on closed subset $F$ by
subsheaf $S_U$ supported on its open complement $U=W-F$.

Technically, this is achieved by introducing new operations-functors
with two inputs and two outputs
\[
\vstretch 80
\begin{tanglec}
\nodeld{n}\Step\noderd{m}\\
\nw1\noder{\co}\node\ne1\\
\nodelu{p}\sw1\se1\noderu{q}
\end{tanglec}
\qquad : \Dbc_{GL(n)\times GL(m)}(pt) \to \Dbc_{GL(p)\times GL(q)}(pt),
\]
which depend on a parameter $\co$ -- a $P_{p,q}\times P_{n,m}$\n-invariant
subset of $GL(n+m)$, where $p+q=n+m$, the parabolic subgroup
$P_{n,m}\subset GL(n+m)$ consists of matrices preserving
$\CC^n\subset\CC^{n+m}$. Minimal such subsets $\co$ are double cosets --
points of the double coset space $P_{p,q}\backslash GL(n+m)/P_{n,m}$.
This is a finite set, it is in bijection with the set of quadruples
$(i,j,k,l)$, $i,j,k,l\in\ZZ_{\ge0}$, which satisfy the equations
$i+j=n$, $k+l=m$, $i+k=p$, $j+l=q$. Hence, we may index the
$P_{p,q}\times P_{n,m}$\n-orbits with these quadruples, say,
$\co_{ijkl}$.

First, we prove that the left and the right hand sides of
\eqref{HnHmHpHq} are isomorphic to the above mentioned operations:
\[
\0\quad
\vstretch 100
\begin{tanglec}
\nodeld{n}\Step\noderd{m}\\
\nw1\nodeld{n+m}\node\ne1\\
\n\\
\nodelu{p}\sw1\se1\noderu{q}
\end{tanglec}
\quad\;\simeq\;\quad
\vstretch 150
\begin{tanglec}
\nodeld{n}\Step\noderd{m}\\
\nw1\noder{GL(n+m)}\node\ne1\\
\nodelu{p}\sw1\se1\noderu{q}
\end{tanglec}
\qquad\qquad;\qquad
\vstretch 60
\begin{tanglec}
\nodeld{n}\step[3]\noderd{m}\\
\nodeld{i\;}\n\noderd{\;j}\step[3]\nodeld{k\;}\n\noderd{\;l}\\
\sw1\se1\step\sw1\se1\\
\id\Step\hxx\Step\id\\
\nw1\ne1\step\nw1\ne1\\
\nodelu{p}\s\step[3]\s\noderu{q}
\end{tanglec}
\quad\;\simeq\;\quad
\vstretch 150
\begin{tanglec}
\nodeld{n}\Step\noderd{m}\\
\nw1\noder{\co_{ijkl}}\node\ne1\\
\nodelu{p}\sw1\se1\noderu{q}
\end{tanglec}
\qquad.
\]
Since $GL(n+m)=\sqcup_{i,j,k,l}\co_{ijkl}$, the former functor above
is a repeated extension of the latter functors via distinguished
triangles.

This shows usefulness of operations with many inputs and outputs for
our purposes. They are also used to prove an equation for associativity
isomorphisms which makes $\cup\Dbc_{GL(k)}(pt)$ into a monoidal
category.  Similar equation for coassociativity isomorphisms is proven
as well.

\begin{acknowledgements}
I am grateful to L.~Crane and D.~Yetter for fruitful discussions which
stimulated this study. Commutative diagrams in this paper are drawn
with the help of a package {\tt diagrams} of Paul Taylor.
\end{acknowledgements}

\numberwithin{equation}{subsection}

\section{Preliminaries}
A definition of equivariant derived categories
is given by Bernstein and Lunts~\cite{BernsL:Equivariant}.
First we explain basic terms.
With a topological space $X$ is associated the category $Sh(X)$
of sheaves of topological spaces.
Its derived category is denoted $D(X)$.
The subcategory consisting of bounded complexes of sheaves
is denoted $D^b(X)$.
If $X$ is a complex algebraic variety,
we call a sheaf \emph{constructible}
if it is constructible with respect to some stratification
by algebraic submanifolds
and stalks are finite-dimensional vector spaces.
A complex is cohomologically constructible
if its cohomology sheaves are constructible.
Subcategory of bounded constructible complexes
is denoted $D^{b,c}(X)$.

\subsection{Equivariant derived categories}
Assume that a complex linear algebraic group $G$
acts algebraically on a complex algebraic variety $X$.
In this setting Bernstein and Lunts~\cite{BernsL:Equivariant}
define bounded constructible equivariant derived category
$D^{b,c}_G(X)$, as a fiber category.

A $G$\n-variety $P$ is called free if $G$ acts freely on $P$
and the quotient map $q:P\to G\backslash P=\overline{P}$
is a locally trivial fibration with the fiber $G$.
A $G$\n-resolution of a $G$\n-variety $X$ is a $G$\n-map $P\to X$,
where the $G$\n-variety $P$ is free.

Let $j:\cj\to\Resto(X,G)$, $P\mapsto (jP:JP\to X)$ be a functor
to the category of $G$\n-resolutions.
Let $\Top$ denote the category of complex algebraic varieties.
Let us denote $\Phi:\Resto(X,G)\to\Top$,
$(R\to X)\mapsto\overline{R}=G\backslash R$ quotient functor.
Consider the composite functor
\[ \Psi:\cj \rTTo^j \Resto(X,G) \rTTo^\Phi \Top, \qquad
P\mapsto (jP:JP\to X) \mapsto \overline{JP}, \]
and define the fiber-category $D^{b,c}(\Psi)$ as follows.

\begin{definition}[Bernstein and Lunts~\cite{BernsL:Equivariant}]
An object of $D^{b,c}(\Psi)$ is a function
$M:\Resto(X,G)\ni P\mapsto M(P)\in D^{b,c}(\overline{JP})$
equipped with isomorphisms
$\alpha_\nu:(\overline{J\nu})^*(M(R))\to M(P)$ given for any
$\nu:P\to R\in\Mor\cj$, such that for any pair of composable morphisms
$P \rTTo^\nu R \rTTo^\mu S$ we have
\[ \alpha_{\mu\nu} = \Bigl[ \overline{J(\mu\nu)}^*M(S) \simeq
\overline{J\nu}^*\overline{J\mu}^*M(S)
\rTTo^{\overline{J\nu}^*\alpha_\mu} \overline{J\nu}^*M(R)
\rTTo^{\alpha_\nu} M(P)\Bigr]. \]
A morphism $\phi:M\to N$ is a collection $\phi(P):M(P)\to N(P)$,
$P\in\Ob\cj$, compatible with $\alpha_\nu$
for any $\nu:P\to R\in\Mor\cj$:
\begin{diagram}
(\overline{J\nu})^*(M(R)) & \rTTo^{\alpha_\nu^M} & M(P) \\
\dTTo<{(\overline{J\nu})^*\phi(R)} & = & \dTTo>{\phi(P)} \\
(\overline{J\nu})^*(N(R)) & \rTTo^{\alpha_\nu^N} & N(P)
\end{diagram}
Define \emph{equivariant derived category} as $\Dbc_G(X)=\Dbc(\Phi)$
in the case of identity functor $j=\id:\cj\rId\Resto(X,G)$.
\end{definition}

We shall also use the notation $\edc{G}X=\Dbc_G(X)$
for equivariant derived category.
Notice that, if $X$ is $G$\n-free, then $\edc{G}X$ is equivalent to
$\Dbc(G\bs X)$.
Without freeness assumption the former and the latter categories
are not equivalent, in general.

Bernstein and Lunts compute the equivariant derived category
in the case when $X$ is a point.

\begin{theorem}[Bernstein and Lunts \cite{BernsL:Equivariant}
Theorem~12.7.2]
Assume that $G$ is a connected linear algebraic group.
The triangulated category $D^{b,c}_G(pt)$
is equivalent to the derived category of the category
of finitely generated differential graded $A$\n-modules,
where the graded algebra $A=H^\bull(BG,\CC)$
is equipped with zero differential.
\end{theorem}

For $G=GL(n,\CC)$ the algebra $A$ is the algebra of symmetric
polynomials of $n$ variables
\[ A_n=\CC[x_1,\dots,x_n]^{\SSS_n} \simeq \CC[e_1,\dots,e_n], \]
where $\deg x_i=2$ and $e_j$ are elementary symmetric functions. For
$G=GL(n_1,\CC)\times\dots\times GL(n_k,\CC)$ we have
$A=A_{n_1}\tens_\CC\dots\tens_\CC A_{n_k}$.

\subsection{Equivariant derived functors}
\subsubsection{The inverse image functor.}
Suppose that $\phi:G\to H$ is a group homomorphism,
and $f:X\to Y$ is a $\phi$\n-equivariant map.
We want to define a functor $f^*=(f,\phi)^*:\Dbc_H(Y)\to\Dbc_G(X)$.
First, denote $\cj=\Resto(f,\phi)$ the category, whose objects
are $\phi$\n-maps $\underline{f}:P\to R$ of resolutions, that is,
\begin{diagram}
P & \rTTo^{\underline{f}} & R \\
\dTTo && \dTTo \\
X & \rTTo^f & Y
\end{diagram}
commutes. Morphisms $\nu:\underline{f}\to\underline{f'}$
are pairs of morphisms of resolutions $\nu_1:P\to P'$, $\nu_2:R\to R'$
such that
\begin{diagram}
P & \rTTo^{\underline{f}} & R \\
\dTTo<{\nu_1} && \dTTo>{\nu_2} \\
P' & \rTTo^{\underline{f'}} & R'
\end{diagram}
commutes. Use the functor $j:\Resto(f;\phi)\to\Resto(X,G)$,
$j(\underline{f}:P\to R)=P$ and $\Psi=j\Phi$.
Bernstein and Lunts \cite{BernsL:Equivariant} have shown that the
restriction functor $\Dbc_g(X)\to\Dbc(\Psi)$ is an equivalence.

Similarly to Bernstein and Lunts \cite{BernsL:Equivariant} we
define the first version of the inverse image functor:
\[ f^*: \Dbc_H(Y) \to \Dbc(\Psi:\Resto(f;\phi)\to\Top), \]
\begin{multline*}
f^*(M:R\mapsto M(R)\in\Dbc(\overline{R})) \\
= \bigl[ f^*M: (P \rTTo^{\underline{f}} R) \mapsto
\overline{\underline{f}}^*(M(R)) \in\Dbc(\overline{P}) \bigr].
\end{multline*}
Next thing is to choose for all $(f,\phi)$ an equivalence
\[ F_{f,\phi}:\Dbc(\Psi:\Resto(f;\phi)\to\Top) \to \Dbc_G(X) \]
quasi-inverse to the canonical restriction functor
\begin{align*}
\Can_f: \Dbc_G(X) &\longrightarrow \Dbc(\Psi:\Resto(f;\phi)\to\Top), \\
[P\mapsto M(P)] &\longmapsto [(f:P\to R)\mapsto M(P)].
\end{align*}
The chosen isomorphisms of the composition
with the identity functor are denoted
%\bigl[ \Dbc(\Resto(f;\phi)\to\Top) \rTTo^{F_{f,\phi}} \Dbc_G(X)
%\rTTo^{\Can_f} \Dbc(\Resto(f;\phi)\to\Top) \bigr] \rTTo^{\eta_f} \Id,
\begin{multline*}
\eta_f:\bigl[ \Dbc(\Resto(f;\phi)\to\Top) \rTTo^{F_{f,\phi}} \Dbc_G(X)
\rTTo^{\Can_f} \Dbc(\Resto(f;\phi)\to\Top) \bigr] \\
\to \Id,
\end{multline*}
\begin{equation*}
\epsilon_f: \Id \to \bigl[ \Dbc_G(X) \rTTo^{\Can_f}
\Dbc(\Resto(f;\phi)\to\Top) \rTTo^{F_{f,\phi}}
\Dbc_G(X) \bigr].
\end{equation*}

Define the \emph{inverse image functor} of $(f,\phi)$ as
\[ f^*=(f,\phi)^*: \Dbc_H(Y) \rTTo^{f^*} \Dbc(\Psi:\Resto(f;\phi)\to\Top)
\rTTo^{F_{f,\phi}} \Dbc_G(X). \]

For composable maps
\[ X \rTTo^f Y \rTTo^g Z \rTTo^h W \quad\text{ over }\quad
G \rTTo^\phi H \rTTo^\psi K \rTTo^\chi L \]
we define categories $\Resto(f,g;\phi,\psi)$ and
$\Resto(f,g,h;\phi,\psi,\chi)$, whose objects are pairs
$(\underline{f},\underline{g})$ (resp. triples
$(\underline{f},\underline{g},\underline{h})$)
of morphisms of resolutions over $(f,g)$ (resp. $(f,g,h)$)
\begin{diagram}
P & \rTTo^{\underline{f}} & R & \rTTo^{\underline{g}} & S &
\rTTo^{\underline{h}} & Q \\
\dTTo && \dTTo && \dTTo && \dTTo \\
X & \rTTo^f & Y & \rTTo^g & Z & \rTTo^h & W
\end{diagram}
Morphisms are triples $P\to P'$, $R\to R'$, $S\to S'$ (resp. quadruples
\dots, $Q\to Q'$) of morphisms of resolutions compatible with
$(\underline{f},\underline{g},\underline{h})$ and
$(\underline{f'},\underline{g'},\underline{h'})$.

There is an isomorphism $\iota:f^*g^*\simeq (gf)^*$ determined uniquely
by the equation
\begin{multline}
\begin{diagram}[inline,nobalance]
&& \Dbc_H(Y) &&&& \Dbc(\Resto(f,g;\phi,\psi)\to\Top) \\
& \ruTTo^{g^*} & \dTwoar<{\iota_{g,f}} & \rdTTo^{f^*} &&& \uTTo>\wr \\
\Dbc_K(Z) && \rTTo_{(gf)^*} && \Dbc_G(X) & \rTTo_\sim &
\Dbc(\Resto(f;\phi)\to\Top)
\end{diagram}
= \\
\begin{diagram}[inline,width=2em,height=1em,objectstyle=\scriptstyle]
\Dbc_H(Y) && \rTTo^{f^*} && \Dbc(\Resto(f;\phi)\to\Top) &&
\rTTo^{F_{f,\phi}} && \Dbc_G(X) \\
& \rdTwoar(2,4)_{\eta_g} \rdTTo(4,4) &&&& \rdId(4,4) &&
\ldTwoar_{\eta_f} &  \\
\uTTo<{F_{g,\psi}} &&&& = && \0 && \dTTo \\
&&&& &&&& \\
\Dbc(\Resto(g;\psi)\to\Top) && \rId && \Dbc(\Resto(g;\psi)\to\Top) &&&&
\Dbc(\Resto(f;\phi)\to\Top) \\
& \rdTwoar(4,4)^{i_{g,f}} &&&& \rdTTo(4,4)^{f^*} &&&  \\
\uTTo<{g^*} &&&& &&&& \dTTo \\
&&&& &&&& \\
\Dbc_K(Z) &&&& \hspace*{-4.5mm} \Dbc(\Resto(gf;\psi\phi)\to\Top)
&& \rTTo_\sim && \Dbc(\Resto(f,g;\phi,\psi)\to\Top) \\
&&& \ruId(4,4) &&&&& \\
\dTTo<{(gf)^*} && \0 && \uTTo>\wr && = && \uTTo>\wr \\
&&& \rdTwoar^{\quad\eta_{gf}^{-1}} &&&&& \\
\Dbc(\Resto(gf;\psi\phi)\to\Top) && \rTTo_{F_{gf,\psi\phi}} && \Dbc_G(X)
&& \rTTo_\sim && \Dbc(\Resto(f;\phi)\to\Top)
\end{diagram}
\label{eq-dia-iota}
\end{multline}

\subsubsection{The direct image functor.}
Let $f:X\to Y$ be a $\phi$\n-map, $\phi:G\to H$.
Assume that $X$ is $K$\n-free, where $K=\Ker\phi\subset G$.
For our purposes it suffices to define $f_*:\Dbc_G(X)\to\Dbc_H(Y)$
as a right adjoint functor to $f^*:\Dbc_H(Y)\to\Dbc_G(X)$.
Furthermore,
we shall use it mainly in the quotient equivalence situation:
$H=G/K$, $\phi:G\to G/K$ is the canonical projection,
$X$ is $K$\n-free, $Y=K\bs X$,
$f=\pi:X\to K\bs X$ is the canonical projection.

\subsubsection{Quotient equivalence.}
\begin{subtheorem}[Bernstein and Lunts~\cite{BernsL:Equivariant}]
Let $K$ be a normal subgroup of $G$,
let $X$ be a $G$\n-space which is free as a $K$\n-space.
Then the quotient map $\pi:X\to K\backslash X$ gives an equivalence
\[ \pi^*: \Dbc_{G/K}(K\backslash X) \to \Dbc_G(X) \]
with a quasi-inverse $\pi_*$.
\end{subtheorem}

In this situation we shall make a concrete choice of a right adjoint
(and quasi-inverse) functor to $\pi^*$
\[ \pi_*: \Dbc_G(X) \to \Dbc_{G/K}(K\bs X), \quad N \mapsto \pi_*N, \]
\begin{multline*}
(\pi_*N)(R\to K\bs X) = N(R\times_{K\bs X}X\to X) \\
\in \Dbc(G\bs(R\times_{K\bs X}X)) \simeq \Dbc((G/K)\bs R),
\end{multline*}
the equivalence is due to the isomorphism
$G\bs(R\times_{K\bs X}X) \simeq (G/K)\bs R$.

\subsubsection{Induction equivalence.}
\begin{subtheorem}[Bernstein and Lunts~\cite{BernsL:Equivariant}]
Let $H$ be a subgroup of $G$, let $X$ be an $H$\n-space.
Then the induction map $i:X\to G\times_HX$, $x\mapsto(1,x)$
gives an equivalence
\[ i^*: D_G(G\times_HX) \to  D_H(X) \]
with a quasi-inverse $i_*$.
\end{subtheorem}

\subsubsection{The direct image with proper supports.}
The following definition belongs to Bernstein and Lunts~
\cite{BernsL:Equivariant}.

\begin{definition}
Let $f:X\to Y$ be a map of $G$\n-varieties.
For any resolution $\pi:P\to Y$ there is a pull-back resolution
$\underline{\pi}$ of $X$.
\begin{diagram}
P\times_YX & \rTTo^{\underline{f}} & P \\
\dTTo<{\underline{\pi}} && \dTTo>{\pi} \\
X & \rTTo^f & Y
\end{diagram}

The functor $f_!:\Dbc_G(X)\to\Dbc_G(Y)$ maps an object
$M:(R\to X)\mapsto M(R\to X)\in\Dbc(\overline{R})$ to the object
$f_!M:(P\to Y)\mapsto
\overline{\underline{f}}_!(M(P\times_YX\to X)) \in\Dbc(\overline{P})$
equipped with an isomorphism
\begin{multline*}
\alpha_\nu^{f_!M}: \overline{\nu}^*\overline{f_R}_!M(R\times_YX\to X)
\rTTo^\beta_\sim
\overline{f_P}_!\overline{\underline{\nu}}^*M(R\times_YX\to X) \to \\
\rTTo^{\overline{f_P}_!\alpha_\nu^M} \overline{f_P}_!M(P\times_YX\to X)
\end{multline*}
for $\nu:P\to R$. Here $\beta$ is a base change isomorphism
obtained from the top square of the following prism.
\begin{diagram}[height=1.5em,width=2.7em,nobalance]
&&&& \overline{R\times_YX} &&& \rTTo^{\overline{f_R}} &&&
\overline{R} \\
&&& \ruTTo(4,2)^{\overline{\underline{\nu}}} & \uTTo &&&&&
\ruTTo(4,2)^{\overline{\nu}} & \\
\overline{P\times_YX} &&& \rTTo^{\overline{f_P}} & \HonV &&
\overline{P} &&&& \uTTo \\
&&&& \vLine && \uTTo &&&& \\
\uTTo &&&& R\times_YX & \hLine & \VonH & \rTTo^{f_R} &&& R \\
&&& \ruTTo(4,2)^{\underline{\nu}} \ldLine(1,2) &&&&&& \ruTTo(4,2)^{\nu}
\ldTTo(2,4) & \\
P\times_YX &&& \HonV & \rTTo^{f_P} && P &&&& \\
& \rdTTo & \ldTTo(1,2) &&&&& \rdTTo &&& \\
&& X &&& \rTTo &&& Y &&
\end{diagram}
Required property of $\alpha$ follows from
\lemref{lem-distr-base-chge} (see Appendix).
\end{definition}

\subsubsection{The equivariant base change isomorphism.}
\label{sec-equi-bas-chg}
Let the pull-back square
\begin{diagram}
W\SEpbk & \rTTo^e & X \\
\dTTo<h && \dTTo>f \\
Z & \rTTo^g & Y
\end{diagram}
consist of equivariant maps:
$h$ is a $G$\n-map, $f$ is an $H$\n-map and $e$, $g$ are $\phi$\n-maps,
where $\phi:G\to H$ is a group homomorphism.
There is a commutative cube
\begin{diagram}[height=1.8em]
P\times_ZW && \rTTo^{\underline{e}} && R\times_YX && \\
& \rdTTo_{\underline{h}} &&& \vLine & \rdTTo^{\underline{f}} & \\
\dTTo && P & \rTTo^{\underline{g}} & \HonV && R \\
&& \dTTo && \dTTo && \\
W & \hLine & \VonH & \rTTo^e & X && \dTTo \\
& \rdTTo_h &&&& \rdTTo^f & \\
&& Z && \rTTo^g && Y
\end{diagram}
where the left and the right walls and the bottom are pull-back squares.
It follows that the top is also a pull-back square.
We define a version of direct image functor with proper support:
\[ h_!: \Dbc(\Resto e\to\Top) \to \Dbc(\Resto g\to\Top),
\quad N\mapsto h_!N, \]
\[ (h_!N)(\underline{g}) = \overline{\underline{h}}_!(N(\underline{e}))
\in\Dbc(\underline{P}), \]
where $N(\underline{e})\in\Dbc(\overline{P\times_ZW})$.
Notice that the functor $h_!$ depends on $f$ as well,
which is not reflected in notation.

The base change isomorphism
\begin{diagram}[width=6em]
\Dbc(\Resto e\to\Top) & \lTTo^{e^*} & \Dbc_H(X) \\
\dTTo<{h_!} & \luTwoar^\beta & \dTTo>{f_!} \\
\Dbc(\Resto g\to\Top) & \lTTo^{g^*} & \Dbc_H(Y)
\end{diagram}
comes from the standard one for quotient spaces. The collection
\[ (g^*f_!M)(\underline{g}) =
\overline{\underline{g}}^*[(f_!M)(R\to Y)] =
\overline{\underline{g}}^*\overline{\underline{f}}_![M(R\times_YX\to X)]
\in\Dbc(\overline{P}) \]
is isomorphic to the collection
\[ (h_!e^*M)(\underline{g}) =
\overline{\underline{h}}_![(e^*M)(\underline{e})] =
\overline{\underline{h}}_!\overline{\underline{e}}^*[M(R\times_YX\to X)]
\in\Dbc(\overline{P}) \]
via base change isomorphism
$\beta:\overline{\underline{g}}^*\overline{\underline{f}}_! \to
\overline{\underline{h}}_!\overline{\underline{e}}^*$.

Finally, we define a full form of base change isomorphism
as the following diagram suggests:
\begin{diagram}[height=2.400000em,width=6em]
\Dbc_G(W) & \pile{\lTTo^{F_{f,\phi}}\\ \Downarrow\eta\\ \rTTo_{\Can_e}}
& \Dbc(\Resto e\to\Top) & \lTTo^{e^*} & \Dbc_H(X) \\
\dTTo<{h_!} & = & \dTTo<{h_!} & \luTwoar^\beta & \dTTo>{f_!} \\
\Dbc_G(Z) &
\pile{\rTTo^{\Can_g}\\ \Downarrow\epsilon \\ \lTTo_{F_{g,\phi}}} &
\Dbc(\Resto g\to\Top) & \lTTo_{g^*} & \Dbc_H(Y)
\end{diagram}
Namely, a full base change isomorphism is:
\begin{multline*}
\beta: g^*_{\text{full}}f_! = F_{g,\phi}g^*f_! \rTTo^{F\beta}
F_{g,\phi}h_!e^* \rTTo^{Fh_!\eta_f^{-1}}
F_{g,\phi}h_!\Can_eF_{f,\phi}e^* \\
= F_{g,\phi}\Can_gh_!F_{f,\phi}e^* \rTTo^{\epsilon_g^{-1}}
h_!F_{f,\phi}e^* = h_!e^*_{\text{full}}.
\end{multline*}

\section{The Hopf category $n_+SL(2)$}
\subsection{Setup and notations}
We partially follow Lusztig \cite{Lus:quiver} and
\cite{Lus:book}, Chapter~9 in notations.
Let $V$ be a vector space and
\[ G = G_V = GL(V) .\]
Let us make the product of $D_G(pt)$ over varying $\dim V$
into a sort of a graded Hopf category.

Assume we are given a decomposition
\[ \cv :\qquad V^1\oplus V^2\oplus\dots\oplus V^k = V \]
into vector subspaces. Associate with it a filtration of $V$
\[ 0=V^{(0)}\subset V^{(1)}\subset \dots\subset V^{(n)}=V ,\qquad
 V^{(m)} = V^1 \oplus\dots\oplus V^m .\]
$P_\cv$ is the corresponding parabolic group
\[ P_\cv = \{ g\in G_V \mid \forall m\ g(V^{(m)}) \subset V^{(m)} \}
= \{ g\in G_V \mid \forall m\ g(V^m) \subset V^{(m)} \} \]
and $U_\cv$ is its unipotent radical. The group
\[ L_\cv = \{ g\in G_V \mid \forall m\ g(V^m) \subset V^m \} =
\prod_{m=1}^k G_{V^m} \]
is a Levi subgroup of $P_\cv$.

Notice that $P_\cv$, $U_\cv$ need only a filtration to be defined
unlike $L_\cv$, which requires a direct sum decomposition.

\subsection{Suggestions for a monoidal 2-category}
To provide a final framework for Hopf categories
one would be interested to have a symmetric monoidal 2\n-category
of equivariant derived categories.
Tensor product of categories would be similar
to that of abelian \kd-linear categories introduced
by Deligne~\cite{Del:cat}. In particular,
$\Dbc_G(X)\dlgn\Dbc_H(Y)\simeq\Dbc_{G\times H}(X\times Y)$
is desirable. However, this wish does not look realistic.
To achieve it, possibly, one has to replace equivariant derived
categories with some other kind of categories, so that inverse
image functors and direct image functors make sense,
and the usual relations still hold.

Let us consider a question of tensor product of functors.
Let $f:X\to Y$, $g:Z\to W$ be maps of algebraic varieties.
Denote by $f^*$, $g^*$, $f_!$, $g_!$
the corresponding equivariant derived functors.
It is explained in \cite{Lyu-sqHopf} that one choice for $f^*\dlgn g^*$
is as good as another as long as they are isomorphic.
If in the isomorphism class of $f^*\dlgn g^*$ there is a functor
$(f\times g)^*$, we can modify the definition of $\dlgn$ and set
$f^*\dlgn g^*=(f\times g)^*$. If $f_!\dlgn g_!$ is isomorphic to
$(f\times g)_!$, we can set $f_!\dlgn g_!=(f\times g)_!$.
The isomorphism
\[ (\Id\dlgn g^*)\circ(f_!\dlgn\Id) \simeq
(f_!\dlgn\Id)\circ(\Id\dlgn g^*) \]
can be chosen as the base change isomorphism
$\beta:(Y\times g)^*(f\times W)_!\to(f\times Z)_!(X\times g)^*$,
constructed for the pull-back square
\begin{diagram}
X\times Z\SEpbk & \rTTo^{X\times g} & X\times W \\
\dTTo<{f\times Z} && \dTTo>{f\times W} \\
Y\times Z & \rTTo^{Y\times g} & Y\times W
\end{diagram}
We stress again that $\dlgn$ is far from being constructed.
Nevertheless, we prove some statements,
which can be interpreted as axioms of a Hopf category.

\subsection{Braiding}
Pretending that the categories
$D_{G_{a_1}\times\dots\times G_{a_k}}(pt)$ form
a monoidal 2\n-category, where $a_i$ are some vector spaces,
we define a braiding in it via functor
\begin{equation*}
c: D_{\prod G_{a_i}\times\prod G_{b_j}}(pt) \rTTo^{[-2d]}
D_{\prod G_{a_i}\times\prod G_{b_j}}(pt)
\rTTo^{\sigma^*} D_{\prod G_{b_j}\times\prod G_{a_i}}(pt) ,
\end{equation*}
where $d=(\sum_i\dim_\CC a_i)(\sum_j\dim_\CC b_j)$ and $\sigma$
is a permutation isomorphism of groups. The shift functor
$K\mapsto K[-2d]$ (the $R$\n-matrix) is related to other functors
that we are using via the following lemma.

\begin{lemma}
Let $h:E\to B$ be an affine linear $C^\infty$\n-bundle.
Then there is a (canonical in $B$) isomorphism of functors
\[ \bigl( \Dbc(B) \rTTo^{h^*} \Dbc(E) \rTTo^{h_!} \Dbc(B) \bigr)
\simeq T^{-2\dim_\CC h}. \]
\end{lemma}

\begin{proof}
Using Propositions 10.1 and 10.8(2) from \cite{Borel-Int_Cohom-V}
we find for any $K\in\Dbc(B)$
\[ h_!h^*K \simeq h_!h^*(\CC\tens K) \simeq h_!(h^*\CC\tens h^*K)
\simeq (h_!h^*\CC)\tens K. \]
Thus we have to prove that $h_!h^*\CC_B\simeq h_!\CC_E$ is isomorphic
to $\CC[-2d]$ for $d=\dim_\CC h$.

Choose a flat connection $\nabla$ on the bundle $h$. By definition
this is a $C^\infty(B)$-linear homomorphism of Lie algebras of
vector fields
\[ \nabla: \Vect(B) \to \Vect(E), \qquad \xi \mapsto \nabla_\xi, \]
such that for each point $e\in E$ we have
$Th({\nabla_\xi}_e)=\xi_{h(e)}$, where $Th:T_eE\to T_{h(e)}B$ is the
tangent map to $h$. The fields $\nabla_\xi$ (the horizontal vector
fields) form a $2b$\n-dimensional (over $\RR$) distribution in $E$,
where $b=\dim_\CC B$. This distribution is involutive, therefore,
by Frobenius theorem, locally there exist coordinates $(z_i,u_j)$,
$i=\overline{1,2b}$, $j=\overline{1,2d}$ in $E$ in which leaves of
the obtained foliation are described by equations $u_1=\const$, \dots,
$u_{2d}=\const$. Moreover, $z_i$ can be chosen as $z_i=x_i\circ h$,
where $x_i$ are local coordinates on the base $B$.

The sheaf of differential graded algebras of \emph{fibrewise forms}
$\Omega^\bull_\nabla$ is defined as a subsheaf of $C^\infty$
differential forms on $E$
\[ \Omega^\bull_\nabla(U) =
\{ \omega\in\Omega^\bull_E(U) \mid \forall \xi\in\Vect(B) \quad
i_{\nabla_\xi}\omega=0 \text{ and } i_{\nabla_\xi}d\omega=0 \}. \]
The second condition may be replaced with $L_{\nabla_\xi}\omega=0$,
since Lie derivative can be computed as
$L_{\nabla_\xi}\omega=di_{\nabla_\xi}\omega+i_{\nabla_\xi}d\omega$.
Therefore, in local coordinates forms $\omega\in\Omega^\bull_\nabla(U)$
are written as $f(u)du^\beta$. Absence of $dz_i$ is implied by the
first condition, and the Lie derivative condition implies independence
of the coefficients on $z_i$ coordinates. We conclude that the complex
$\Omega^\bull_\nabla$ is a $c$\n-soft on fibers of $h$ resolution of
the constant sheaf $\CC_E$. Hence, it can be used to compute
the complex $h_!\CC_E$.

Let $V$ be an open subset of $B$. Then $h_!\CC_E(V)$ is a complex of
vector spaces $\Omega^\bull_{\nabla,c}(h^{-1}(V))$, where $c$ indicates
such forms $\omega$ that $h^{-1}(K)\cap\supp\omega$ is compact for each
compact subset $K$ of $B$. This complex is exact everywhere except the
maximal degree $2d$. It has a map into the algebra of functions on $V$
\[ \alpha(V): \Omega^{2d}_{\nabla,c}(h^{-1}(V)) \to C^\infty(V), \quad
\omega \longmapsto (v\mapsto \int_{h^{-1}(v)}\omega), \]
given by fibrewise integration. Local presentation of $\omega$ implies
that the function $\alpha(V)(\omega)$ is locally constant, hence,
is in $\CC_B(V)$. The obtained chain map $\alpha:h_!\CC_E\to\CC_B[-2d]$
is a quasi-isomorphism.
\end{proof}

Another (complex analytic) construction of $\alpha$ will be published
elsewhere. It follows the lines of \cite{Lyu-l-adic-Mat-Studii}, where
the case of quasicompact schemes over $\overline{\FF_p}$ and
$\ell$\n-adic sheaves is considered. In the analytic setting one can
also prove that for any pull-back square
\begin{diagram}
F \SEpbk & \rTTo^j & A \\
\dTTo<g && \dTTo>f \\
E & \rTTo^h & B
\end{diagram}
where $h$ is an affine linear bundle (and so is $j$), we have an equation
\[
\begin{diagram}[inline,width=2.5em]
&& \Dbc(E) && \\
& \ruTTo^{h^*} & \dTTo>{g^*} & \rdTTo^{h_!} & \\
\Dbc(B) && \Dbc(F) && \Dbc(B) \\
\dTTo<{f^*} & \ruTTo^{j^*} & \dTwoar>{\alpha_A}
& \rdTTo^{j_!} & \dTTo>{f^*} \\
\Dbc(A) && \rTTo_{[-2d]} && \Dbc(A)
\end{diagram}
\quad = \quad
\begin{diagram}[inline,width=2.5em]
&& \Dbc(E) && \\
& \ruTTo^{h^*} & \dTwoar>{\alpha_B} & \rdTTo^{h_!} & \\
\Dbc(B) && \rTTo_{[-2d]} && \Dbc(B) \\
\dTTo<{f^*} && = && \dTTo>{f^*} \\
\Dbc(A) && \rTTo_{[-2d]} && \Dbc(A)
\end{diagram}
\]
This is the precise meaning of canonicity of isomorphism $\alpha$.

\begin{corollary}\label{cor-aff-bun-equi}
Let $h:E\to B$ be an affine linear $G$\n-bundle (an affine linear bundle
equipped with a group homomorphism $G\to\Aut_{\text{af.lin.bun.}}(h)$).
%a $G$\n-equivaraint map, which is an affine linear bundle, such that the
%action of any element $g\in G$ is an automorphism of affine linear bundle $h$).
%linear combination maps $E\times_BE\to E$,
%$(e,f)\mapsto \lambda e+(1-\lambda)f$ are
%$G$\n-equivariant for all $\lambda\in\CC$).
Then there is an isomorphism of functors
\[ \bigl( \Dbc_G(B) \rTTo^{h^*} \Dbc_G(E) \rTTo^{h_!} \Dbc_G(B) \bigr)
\simeq T^{-2\dim_\CC h}. \]
\end{corollary}

\begin{proof}
The system of isomorphisms $\alpha$ in $\Dbc(\overline{P})$ for
$G$\n-resolutions $P\to X$ is coherent due to canonicity of $\alpha$.
\end{proof}

\subsection{Operations}
Let two decompositions of $V$ into a direct sum be given:
\begin{alignat}2
\cv &:&\qquad V^1\oplus V^2\oplus\dots\oplus V^k &= V
\label{eq-decomp-V-to-Vi2} \\
\cw &:&\qquad W^1\oplus W^2\oplus\dots\oplus W^l &= V.
\label{eq-decomp-V-to-Wj2}
\end{alignat}
Let $\co\subset G$ be a left $P_\cw$\n-invariant and right
$P_\cv$\n-invariant subset. We associate with it an operation
\[ \zhe_\cw^{\co;\cv} = \quad
\begin{tangle}
\vstretch 120
\object{V^1}\Step\object{V^2}\step[6]\object{V^k} \\
\nw4\nw0\nw2\nodel{\co}\nw0\n\ne0\ne2\ne0\ne4 \\
\sw4\sw0\sw2\sw0\id\se0\se2\se0\se4 \\
\object{W^1}\Step\object{W^2}\step[6]\object{W^k}
\end{tangle}
\quad = \quad
\vstretch 68
\begin{tangle}
\object{V^1}\Step\object{V^2}\step[6]\object{V^k} \\
\nw4\nw0\nw2\nodel{\co}\nw0\n\ne0\ne2\ne0\ne4 \\
\step[4]\idash \\
\sw4\sw0\sw2\sw0\s\se0\se2\se0\se4 \\
\object{W^1}\Step\object{W^2}\step[6]\object{W^k}
\end{tangle}
\quad = \pitchfork_\cw^\cv \circ \Psi_\cw^\cv.
\]
The components of it are defined below.

\begin{notations}
In graphical notations a straight line $\Big|X$ labeled by a vector
space $X$ denotes category $D_{G_X}(pt)$. A dashed line
$
\begin{tangle}
\idash
\end{tangle}
\:\cx$ labeled by a filtration $\cx$ of a vector space $X$ denotes
category $D_{P_\cx}(pt)$. If a horizontal line crosses a diagram
intersecting transversally with several (solid or dashed) open intervals
labeled by spaces $Y_1$, $Y_2$,\dots and filtrations $\cx_1$, $\cx_2$
of spaces $X_1$, $X_2$,\dots , this line represents the category
\( D_{G_{Y_1}\times G_{Y_2}\times\dots\times
P_{\cx_1}\times P_{\cx_2}\times\dotsb}(pt) \).
The order of factors repeats the order of intersection points on the
horizontal line.

A vertex (or a horizontal row of vertices) represents a functor from
the category, corresponding to horizontal line just above the vertex,
to the category represented by a line just below the vertex.
\end{notations}

\subsubsection{The multiplication.}\label{sec-Multiplication}
Multiplication operation is
\begin{multline*}
\qquad
\vstretch 60
\begin{tangle}
\nodelu{V^1}\Step\nodelu{V^2}\step[6]\noderu{V^k} \\
\nw4\nw0\nw2\nodel{\co}\nw0\n\ne0\ne2\ne0\ne4 \\
\step[4]\idash\noder{\cw} \\
\end{tangle}
\quad = \Psi_\cw^{\co;\cv} = \Bigl( D_{L_\cv}(pt) \rTTo^{\phi^*}
D_{P_\cw\times L_\cv}(\co/{U_\cv}) \\
\rTTo^{\pi_*} D_{P_\cw}(\co/{P_\cv}) \rTTo^{\alpha_!}
D_{P_\cw}(pt) \Bigr),
\end{multline*}
Here the scheme of multiplication is similar to that in \cite{Lus:book}:
\begin{equation}
pt \lTTo^\phi \co/{U_\cv} \rTTo^\pi \co/{P_\cv} \rTTo^\alpha pt,
\label{eq-p1p2p32}
\end{equation}
where $\pi$ is a canonical projection.
The action of $L_\cv=P_\cv/U_\cv$ in $\co/{U_\cv}$ is
induced from the action $p.o=op^{-1}$, $p\in P_\cv$.

In particular, for $l=1$ we have $P_\cw=G_V$, $\co=G_V$ and
the multiplication operation is
\begin{multline*}
\qquad
\vstretch 60
\begin{tangle}
\nodelu{V^1}\Step\nodelu{V^2}\step[6]\noderu{V^k} \\
\nw4\nw0\nw2\nodel{\co}\nw0\n\ne0\ne2\ne0\ne4 \\
\step[4]\id \\
\end{tangle}
\quad = \Psi_1^{\cv} = \Bigl( D_{L_\cv}(pt) \rTTo^{\phi^*}
D_{G_V\times L_\cv}(G_V/{U_\cv}) \\
\rTTo^{\pi_*} D_{G_V}(G_V/{P_\cv}) \rTTo^{\alpha_!} D_{G_V}(pt) \Bigr),
\end{multline*}
Here the scheme of multiplication is precisely that of Lusztig~
\cite{Lus:book}:
\begin{equation*}
pt \lTTo^\phi G_V/{U_\cv} \rTTo^\pi G_V/{P_\cv} \rTTo^\alpha pt.
\end{equation*}

The particular case $k=1$, $\co=G_V$ is also important. We have then
$L_\cv=P_\cv=G_V$, $U_\cv=1$, and
\[ \Psi^1_\cw = \bigl( D_G(pt) \rTTo^{\pr_2^*}
D_{P_\cw\times G}(G) \rTTo^{d_{1*}} D_{P_\cw}(pt) \bigr) ,\]
where $d_1:G\to pt$.

\begin{proposition}\label{pro-psi1=res}
The functor $\Psi^1_\cw$ is isomorphic to the restriction functor
$\Resto_{P_\cw,G}:D_G(pt)\to D_{P_\cw}(pt)$.
\end{proposition}

\begin{proof}
A map $d_1$ is the map $G\rEpi pt$.
Hence, $d_{1*}: D_{P_\cw\times G}(G) \to D_{P_\cw}(pt)$
is an equivalence with a quasi-inverse
$d_1^*: D_{P_\cw}(pt) \to D_{P_\cw\times G}(G)$.

The map $s_0:pt\to G$, $s_0(pt)=1$, is a $\Delta$\n-map
with respect to the homomorphism of groups
$\Delta:P_\cw\to P_\cw\times G$, $\Delta(p)=(p,p)$.
Furthermore, $s_0$ is an induction map
$pt \rTTo^i (P\times G)/P \rTTo^\rho_\sim G$, where $i(pt)=(1,1)$,
$\rho(p,g)=pg^{-1}$ and $\rho$ is a $P\times G$\n-map.
Therefore, $s_0^*: D_{P_\cw\times G}(G) \to D_{P_\cw}(pt)$
is an induction equivalence. Since
\( \bigl( pt \rTTo^{s_0} G \rTTo^{d_1} pt \bigr) = \id\),
we have
\[ \bigl( D_{P_\cw}(pt) \rTTo^{d_1^*} D_{P_\cw\times G}(G)
\rTTo^{s_0^*} D_{P_\cw}(pt) \bigr) \simeq \Id .\]
Therefore, $s_0^* \simeq d_{1*}$. Hence,
\[ \Psi^1_\cw \simeq \bigl( D_G(pt) \rTTo^{\pr_2^*}
D_{P_\cw\times G}(G) \rTTo^{s_0^*} D_{P_\cw}(pt) \bigr)
\simeq \Resto_{P_\cw,G},\]
for $\bigl(pt \rTTo^{s_0} G \rTTo pt\bigr) = \id$ and
$P_\cw \rTTo^\Delta P_\cw\times G \rTTo^{\pr_2} G$ is an inclusion.
\end{proof}

\subsection{Associativity}
Assume that besides decomposition \eqref{eq-decomp-V-to-Vi2} of $V$ into
direct sum of $V^i$ we have also decompositions
\[ \cv^i:\qquad V^i = \bigoplus_{m=1}^{m^i} V_m^i .\]
We can produce out of these the decomposition
\[ \cu:\qquad V = \oplus_{i,m} V_m^i ,\]
where the lexicographic order of summands is used
$V_1^1$, \dots,$V_{m^1}^1,V_1^2$, \dots,$V_{m^2}^2$,
\dots,$V_1^l$, \dots,$V_{m^l}^l$.

An associativity isomorphism is obtained from the isomorphism
\[
%\Psi^{\cv^1}\dlgn\Psi^{\cv^2}\dlgn\dots\dlgn\Psi^{\cv^k}\circ\Psi^\cv = \quad
\vstretch 80
\begin{tangle}
\nw1\nodelu{\cv^1}\n\ne1\step\nw1\nodelu{\cv^2}\n\ne1\step%
\nw1\n\noderu{\cv^k}\ne1 \\
\step\nw3\step\nodel{\co}\step\n\step\noder{\cv}\step\ne3 \\
\step[4]\idash
\end{tangle}
\quad \rTTo^\sim \quad
\vstretch 120
\begin{tangle}
\nw4\nw3\nw2\nw1\n\ne1\noderu{\mathcal U}\ne2\ne3\ne4 \\
\step[4]\idash
\end{tangle}
\quad = \Psi^\cu.
\]
This isomorphism can be read from the diagram in \figref{fig-ass}.

\begin{figure}
\begin{diagram}[height=2.400000em]
\eDc{\prod L_{\cv^i}}{pt} && \rId && \eDc{L_\cu}{pt} \\
\dTTo<{\prod\phi^{i*}} &&&& \dTTo>{\phi_\cu^*} \\
\eDc{\prod{G_{V^i}\times L_{\cv^i}}}{\prod_iG_{V^i}/{U_{\cv^i}}}
& \rTTo^{\pr_{Y2}^*} &
\eDc{P_\cw\times L_\cv\times L_\cu}
{\co/{U_\cv}\times\prod_iG_{V^i}/{U_{\cv^i}}}
& \rTTo^{\pi'_*}_\sim & \eDc{P_\cw\times L_\cu}{\co/{U_\cu}} \\
\dTTo<{\prod\pi_*^i}>\wr && \dTTo<{\pi''_*}>\wr &&
\dTTo<\wr>{\pi_{\cu*}} \\
\eDc{\prod G_{V^i}}{\prod_iG_{V^i}/{P_{\cv^i}}} &
\rTTo^{\pr_{X2}^*} &
\eDc{P_\cw\times L_\cv}{\co/{U_\cv}\times\prod_iG_{V^i}/{P_{\cv^i}}}
& \rTTo^{\pi'''_*}_\sim & \eDc{P_\cw}{\co/{P_\cu}} \\
\dTTo<{\prod\alpha^i_!} && \dTTo<{\pr_{X1!}} && \dTTo>{\beta_!} \\
\eDc{L_\cv}{pt} & \rTTo^{\phi_\cv^*} & \eDc{P_\cw\times L_\cv}{\co/{U_\cv}}
& \rTTo^{\pi_{\cv*}}_\sim & \eDc{P_\cw}{\co/{P_\cv}} \\
&&&& \dTTo>{\alpha_{\cv!}} \\
&&&& \eDc{P_\cw}{pt}
\end{diagram}
\caption{Associativity isomorphism}
\label{fig-ass}
\end{figure}
%% end try

The action of $P_\cw\times L_\cv$ in
$X=\co/{U_\cv}\times\prod_iG_{V^i}/{P_{\cv^i}}$ is the following
\[ (p,l).[o\times\prod g_i] =
pol^{-1} \times\prod l_ig_i, \quad l=\prod l_i .\]
The map $\pi'''$ is the quotient map
\[ \pi''': X \to L_\cv\backslash X \simeq \co/{P_\cu} .\]
The action of $P_\cw\times L_\cv\times L_\cu$ in
$Y=\co/{U_\cv}\times\prod_iG_{V^i}/{U_{\cv^i}}$ is the following
\[ (p,l,\bar l).[o \times\prod g_i] =
pol^{-1} \times\prod l_ig_i\bar l_i^{-1}, \quad
l=\prod l_i,\;\bar l=\prod\bar l_i .\]
The maps $\pi'$, $\pi''$ are the quotient maps
\[ \pi': Y \to L_\cv\backslash Y \simeq \co/{U_\cu} ,\quad
\pi'': Y \to L_\cu\backslash Y \simeq X .\]
The canonical projection $\beta$ comes from the inclusion
$P_\cu\subset P_\cv$,
\[ \beta: \co/{P_\cu} \to \co/{P_\cv} .\]

\subsection{The associativity equation}
Assume that we have the following decompositions
of finite-dimensional $\CC$\n-vector spaces:
\begin{alignat*}4
\cw &:&\quad V &= \oplus_j W^j, &\qquad
\cv &:&\quad V &= \oplus_i V^i, \\
\cv^i &:&\quad V^i &= \oplus_m V^i_m, &\qquad
\cv^i_m &:&\quad V^i_m &= \oplus_p V^i_{m,p}.
\end{alignat*}
These decompositions imply also the following decompositions:
\begin{equation*}
\cu :\; V = \oplus_{i,m} V^i_m, \quad\;
\cy^i :\; V^i = \oplus_{m,p} V^i_{m,p}, \quad\;
\cx :\; V = \oplus_{i,m,p} V^i_{m,p}.
\end{equation*}
Let $\co\subset G_V$ be a left $P_\cw$\n-invariant
and right $P_\cv$\n-invariant subset.
Associativity isomorphisms give two isomorphisms
between the composite operation and the single operation
as in diagram
\begin{diagram}[height=4.5em,LaTeXeqno]
\vstretch 75
\begin{tangle}
\nw1\nodel{\cv^1_1\;}\n\step\nw1\n\step\nw1\n\step\nw1\n\ne1 \\
\step\nw1\nodel{\cv^1\;}\node\ne1\Step\nw1\node\ne1 \\
\Step\nw2\step\nodelu{\cv\;}\nodeld{\cw\;}\noder{\;\co}\step\ne2 \\
\Step\sw2\sw1\s\se1\se2
\end{tangle}
& \rTTo^{\assoc} &
\begin{tangle}
\nw2\nw1\nodel{\cy^1\;}\n\ne1\step\nw2\nw1\n\ne1\ne2 \\
\Step\nw2\step\nodelu{\cv\;}\nodeld{\cw\;}\noder{\;\co}\step\ne2 \\
\Step\sw2\sw1\s\se1\se2
\end{tangle}
\\
\dTTo<{\assoc} & = & \dTTo>{\assoc} & . \\
\begin{tangle}
\nw1\nodel{\cv^1_1\;}\n\step\nw1\n\step\nw1\n\step\nw1\n\ne1 \\
\step\nw3\nw0\nw1\nodelu{\CU\;\;\;}\nodeld{\cw\;}\noder{\;\co}
\ne1\ne0\ne3\\
\Step\sw2\sw1\s\se1\se2
\end{tangle}
& \rTTo^{\assoc} &
\vstretch 150
\begin{tangle}
\nw4\nw3\nw2\nw1\nodelu{\cx\;\;\;}\nodeld{\cw\;}\noder{\;\co}
\n\ne1\ne2\ne3\ne4\\
\Step\sw2\sw1\id\se1\se2
\end{tangle}
\label{dia-asso-square}
\end{diagram}

\begin{proposition}[Associativity equation]
Diagram \eqref{dia-asso-square} is commutative.
\end{proposition}

\begin{proof}
The left-lower associativity composition is given on diagram in
\figref{fig-ll-assoc-start}.
\begin{figure}
\begin{diagram}[width=3.5em,height=2.400000em,nobalance]
\edc{\prod G_{V^i_{mp}}}{pt} & \rTTo^{\phi_\cx^*} &
\edc{P_\cw\times L_\cx}{\co/U_\cx} & \rTTo^{\pi_{\cx*}}_\sim &
\edc{P_\cw}{\co/P_\cx} & \rTTo^{\alpha_{\cx!}} & \edc{P_\cw}{pt} \\
\dTTo<{\prod\phi_m^{i*}} && \uTTo<\sim>{\pi'_*} &
\ruTTo(2,4)^\sim_{\pi''_*} && \ruTTo(2,4)_{\alpha_{\cu!}} & \\
\edc{\prod G_{V^i_m}L_{\cv^i_m}}{\prod G_{V^i_m}/U_{\cv^i_m}} &
\rTTo^{\pr_2^*} &
\edc{P_\cw\times L_\cu\times L_\cx}
{\co/U_\cu\times\prod G_{V^i_m}/U_{\cv^i_m}}
&& \dTTo>{\beta_{\cx!}} && \uTTo>{\alpha_{\cv!}} \\
\dTTo<{\prod\pi^I_{m*}}>\wr && \dTTo<\wr>{\pi''_*} &&&& \\
\edc{\prod G_{V^i_m}}{\prod G_{V^i_m}/P_{\cv^i_m}} & \rTTo^{\pr_2^*} &
\edc{P_\cw\times L_\cu}{\co/U_\cu\times\prod G_{V^i_m}/P_{\cv^i_m}} &&
\edc{P_\cw}{\co/P_\cu} & \rTTo^{\beta_!} & \edc{P_\cw}{\co/P_\cv} \\
&& \dTTo>{\pr_{1!}} & \ruTTo^{\pi_{\cu*}}_\sim &&& \\
&& \edc{P_\cw\times L_\cu}{\co/U_\cu} && \uTTo<\wr>{\pi'''_*} &&
\uTTo<\wr>{\pi_{\cv*}} \\
\dTTo<{\prod\alpha^i_{m!}} & \ruTTo(2,4)^{\phi^*_\cu} &
\uTTo<\wr>{\pi'_*} &&&& \\
&& \edc{P_\cw\times L_\cv\times L_\cu}
{\co/U_\cv\times\prod G_{V^i}/U_{\cv^i}} &
\rTTo^{\pi''_*}_\sim &
\edc{P_\cw\times L_\cv}{\co/U_\cv\times\prod G_{V^i}/P_{\cv^i}} &
\rTTo^{\pr_{1!}} & \edc{P_\cw L_\cv}{\co/U_\cv} \\
&& \uTTo>{\pr_2^*} && \uTTo>{\pr_2^*} && \uTTo>{\phi_\cv^*} \\
\edc{\prod G_{V^i_m}}{pt} & \rTTo^{\nquad\prod\phi^{i*}\nquad} &
\edc{\prod G_{V^i}\times L_{\cv^i}}{\prod G_{V^i}/U_{\cv^i}} &
\rTTo^{\nquad\prod\pi_*^i\nquad}_\sim &
\edc{\prod G_{V^i}}{\prod G_{V^i}/P_{\cv^i}} &
\rTTo^{\nquad\prod\alpha_{i!}\nquad} & \edc{L_\cv}{pt}
\end{diagram}
\caption{First composition of associativity isomorphisms}
\label{fig-ll-assoc-start}
\end{figure}
The right-upper associativity composition is given on diagram in
\figref{fig-ru-assoc-start}. We have to check that isomorphism in this
figure equals to the one in \figref{fig-ll-assoc-start}.

\begin{figure}
\begin{diagram}[width=2.6em,height=2.400000em,inline]
\edc{\prod G_{V^i_{mp}}}{pt} & \rTTo^{\phi_\cx^*} &
\edc{P_\cw\times L_\cx}{\co/U_\cx} & \rTTo^{\pi_{\cx*}}_\sim &
\edc{P_\cw}{\co/P_\cx} & \rTTo^{\alpha_{\cx!}} & \edc{P_\cw}{pt} \\
& \rdTTo(2,4)_{\prod\phi^*_{\cy^i}} & \uTTo<\wr>{\pi'_*} &&&
\rdTTo^{\beta_!} & \uTTo>{\alpha_{\cv!}} \\
&& \edc{P_\cw\times L_\cv\times L_\cx}{\co/U_\cv\times\prod
G_{V^i}/U_{\cv^i}} && \uTTo<\wr>{\pi'''_*} && \edc{P_\cw}{\co/P_\cv} \\
\dTTo<{\prod\phi_m^{i*}} && \uTTo>{\pr_2^*} & \rdTTo_\sim^{\pi''_*} &&&
\uTTo<\wr>{\pi_{\cv*}} \\
&& \edc{\prod G_{V^i}\times L_{\cy^i}}{\prod G_{V^i}/U_{\cy^i}} &&
\hspace*{-3mm} \edc{P_\cw\times L_\cv}
{\co/U_\cv\times\prod G_{V^i}/P_{\cy^i}} &
\rTTo^{\pr_{1!}} & \edc{P_\cw\times L_\cv}{\co/U_\cv} \\
&& \uTTo<\wr>{\pi'_*} & \rdTTo(2,4)_\sim^{\prod\pi_{\cy^i*}} &&& \\
\edc{\prod G_{V^i_m}\times L_{\cv^i_m}}{\prod G_{V^i_m}/U_{\cv^i_m}} &
\rTTo^{\pr_2^*} & \edc{\prod G_{V^i}\times L_{\cv^i}\times
L_{\cy^i}}{\prod G_{V^i}/U_{\cv^i}\times\prod G_{V^i_m}/U_{\cv^i_m}} &&
\uTTo>{\pr_2^*} && \\
\dTTo<{\prod\pi^i_{m*}}>\wr && \dTTo<\wr>{\pi''_*} &&&&
\uTTo>{\phi_\cv^*} \\
\edc{\prod G_{V^i_m}}{\prod G_{V^i_m}/P_{\cv^i_m}} & \rTTo^{\pr_2^*} &
\edc{\prod G_{V^i}\times L_{\cv^i}}{\prod G_{V^i}/U_{\cv^i}\times\prod
G_{V^i_m}/P_{\cv^i_m}} & \rTTo^{\pi'''_*}_\sim &
\edc{\prod G_{V^i}}{\prod G_{V^i}/P_{\cy^i}} && \\
\dTTo<{\prod\alpha^i_{m!}} && \dTTo>{\pr_{1!}} &&
\dTTo>{\prod\beta_{i!}} & \rdTTo^{\prod\alpha_{\cy^i!}} & \\
\edc{\prod G_{V^i_m}}{pt} & \rTTo^{\nquad\prod\phi^{i*}\nquad} &
\edc{\prod G_{V^i}\times L_{\cv^i}}{\prod G_{V^i}/U_{\cv^i}} &
\rTTo^{\nquad\prod\pi_*^i\nquad}_\sim &
\edc{\prod G_{V^i}}{\prod G_{V^i}/P_{\cv^i}} &
\rTTo^{\nquad\prod\alpha_{i!}\nquad} & \edc{L_\cv}{pt}
\end{diagram}
\caption{Second composition of associativity isomorphisms}
\label{fig-ru-assoc-start}
\end{figure}

The subdiagram of diagram in \figref{fig-ru-assoc-start}
between the third and the forth columns is transformed via
\propref{prop-distr-base-chge} to a 3\n-column subdiagram.
The right part of it coincides with the right subdiagram in
\figref{fig-ll-assoc-start} between the third and the fourth columns
and can be canceled.
The left subdiagram of \figref{fig-ll-assoc-start} between the first
and the third columns is also transformed via
\propref{prop-distr-base-chge}.
We come to equation between the isomorphisms in Figures
\ref{fig-ll-assoc-1} and \ref{fig-ru-assoc-1}.

\begin{figure}
\rotatebox{90}{%
\hspace*{1mm}%
\begin{diagram}[width=4.0em,height=2.600000em,inline]
&&&&&& \\
&&&&&& \\
\edc{\prod G_{V^i_{mp}}}{pt} & \rTTo^{\phi_\cx^*} &
\edc{P_\cw\times L_\cx}{\co/U_\cx} && \rTTo^{\pi_{\cx*}}_\sim &&
\edc{P_\cw}{\co/P_\cx} \\
\dTTo<{\prod\phi_m^{i*}} && \uTTo<\sim>{\pi'_*} &&&
\ruTTo_\sim^{\pi'''_*} \ruTTo(2,4)_{\pi_*}^\sim & \\
\edc{\prod G_{V^i_m}\times L_{\cv^i_m}}{\prod G_{V^i_m}/U_{\cv^i_m}} &
\rTTo^{\pr_2^*} &
\edc{P_\cw\times L_\cu\times L_\cx}{\co/U_\cu\times\prod
G_{V^i_m}/U_{\cv^i_m}} & \rTTo^{\pi''_*}_\sim &
\edc{P_\cw\times L_\cu}{\co/U_\cu\times\prod G_{V^i_m}/P_{\cv^i_m}} &&
\dTTo>{\beta_{\cx!}} \\
&& \dTTo<\wr>{\pi''_*} & \ruId & \uTTo<\wr>{\pi'_*} && \\
\dTTo<{\prod\pi^i_{m*}}>\wr &&
\edc{P_\cw\times L_\cu}{\co/U_\cu\times\prod G_{V^i_m}/P_{\cv^i_m}} &
\rTTo_\sim^{\pi^{\prime*}} &
\edc{P_\cw\times L_\cv\times L_\cu}{\co/U_\cv\times\prod
G_{V^i}/U_{\cv^i}\times\prod G_{V^i_m}/P_{\cv^i_m}} &&
\edc{P_\cw}{\co/P_\cv} \\
& \ruTTo^{\pr_2^*} && \ruTTo(4,2)_{\pr_3^*} & \dTTo>{\pr_{12!}} &
\ruTTo^{\pi_*}_\sim & \uTTo<\wr>{\pi'''_*} \\
\edc{\prod G_{V^i_m}}{\prod G_{V^i_m}/P_{\cv^i_m}} &&&&
\edc{P_\cw\times L_\cv\times L_\cu}{\co/U_\cv\times\prod
G_{V^i}/U_{\cv^i}} & \rTTo^{\pi''_*}_\sim &
\edc{P_\cw\times L_\cv}{\co/U_\cv\times\prod G_{V^i}/P_{\cv^i}} \\
\dTTo<{\prod\alpha^i_{m!}} &&& \ruTTo(4,2)^{\phi^*} & \uTTo>{\pr_2^*}
&& \uTTo>{\pr_2^*} \\
\edc{\prod G_{V^i_m}}{pt} && \rTTo^{\prod\phi^{i*}} &&
\edc{\prod G_{V^i}\times L_{\cv^i}}{\prod G_{V^i}/U_{\cv^i}} &
\rTTo^{\nquad\prod\pi_*^i\nquad}_\sim &
\edc{\prod G_{V^i}}{\prod G_{V^i}/P_{\cv^i}}
\end{diagram}
\hspace*{1mm}
}
\caption{First isomorphism}
\label{fig-ll-assoc-1}
\end{figure}

\begin{figure}
\rotatebox{90}{%
\hspace*{1mm}%
\begin{diagram}[width=4.0em,height=2.4500000em,inline]
&&&&&& \\
\edc{\prod G_{V^i_{mp}}}{pt} & \rTTo^{\phi_\cx^*} &
\edc{P_\cw\times L_\cx}{\co/U_\cx} & \rTTo^{\pi_{\cx*}}_\sim &
\edc{P_\cw}{\co/P_\cx} & \rTTo^{\beta_{\cx!}} & \edc{P_\cw}{\co/P_\cu}
\\
& \rdTTo(2,4)_{\prod\phi^*_{\cy^i}} & \uTTo<\wr>{\pi'_*} &&
\uTTo<\wr>{\pi'''_*} && \uTTo<\wr>{\pi'''_*} \\
&& \edc{P_\cw\times L_\cv\times L_\cx}{\co/U_\cv\times\prod
G_{V^i}/U_{\cy^i}} & \rTTo^{\pi''_*} &
\edc{P_\cw\times L_\cv}{\co/U_\cv\times\prod G_{V^i}/P_{\cy^i}} &
\rTTo^{\beta_!} &
\edc{P_\cw\times L_\cv}{\co/U_\cv\times\prod G_{V^i}/P_{\cv^i}} \\
\dTTo<{\prod\phi_m^{i*}} && \uTTo>{\pr_2^*} &&&& \\
&& \edc{\prod G_{V^i}\times L_{\cy^i}}{\prod G_{V^i}/U_{\cy^i}} &&&& \\
&& \uTTo<\wr>{\pi'_*} & \rdTTo(2,4)_\sim^{\prod\pi_{\cy^i*}} &
\uTTo>{\pr_2^*} && \\
\edc{\prod G_{V^i_m}\times L_{\cv^i_m}}{\prod G_{V^i_m}/U_{\cv^i_m}} &
\rTTo^{\pr_2^*} & \edc{\prod G_{V^i}\times L_{\cv^i}\times
L_{\cy^i}}{\prod G_{V^i}/U_{\cv^i}\times\prod G_{V^i_m}/U_{\cv^i_m}} &&
&& \uTTo>{\pr_2^*} \\
\dTTo<{\prod\pi^i_{m*}}>\wr && \dTTo<\wr>{\pi''_*} &&&& \\
\edc{\prod G_{V^i_m}}{\prod G_{V^i_m}/P_{\cv^i_m}} & \rTTo^{\pr_2^*} &
\edc{\prod G_{V^i}\times L_{\cv^i}}{\prod G_{V^i}/U_{\cv^i}\times\prod
G_{V^i_m}/P_{\cv^i_m}} & \rTTo^{\pi'''_*}_\sim &
\edc{\prod G_{V^i}}{\prod G_{V^i}/P_{\cy^i}} && \\
\dTTo<{\prod\alpha^i_{m!}} && \dTTo>{\pr_{1!}} &&&
\rdTTo>{\prod\beta_{i!}} & \\
\edc{\prod G_{V^i_m}}{pt} & \rTTo^{\nquad\prod\phi^{i*}\nquad} &
\edc{\prod G_{V^i}\times L_{\cv^i}}{\prod G_{V^i}/U_{\cv^i}} &&
\rTTo^{\prod\pi_*^i}_\sim &&
\edc{\prod G_{V^i}}{\prod G_{V^i}/P_{\cv^i}}
\end{diagram}
\hspace*{1mm}
}
\caption{Second isomorphism}
\label{fig-ru-assoc-1}
\end{figure}

\figref{fig-ll-assoc-1} can be transformed further using
\propref{prop-distr-base-chge} to the form of \figref{fig-ll-assoc-2}.
We have to prove that it equals to isomorphism in
\figref{fig-ru-assoc-1}.

\begin{figure}
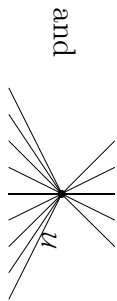

\rotatebox{90}{%
\hspace*{1mm}%
\begin{diagram}[width=4.0em,height=2.600000em,inline]
&&&&&& \\
\edc{\prod G_{V^i_{mp}}}{pt} & \rTTo^{\phi_\cx^*} &
\edc{P_\cw\times L_\cx}{\co/U_\cx} &&&& \\
\dTTo<{\prod\phi_m^{i*}} && \uTTo<\sim>{\pi'_*}
& \rdTTo(4,2)^{\pi_{\cx*}}_\sim &&& \\
\edc{\prod G_{V^i_m}\times L_{\cv^i_m}}{\prod G_{V^i_m}/U_{\cv^i_m}}
& \rTTo^{\pr_2^*} &
\edc{P_\cw\times L_\cu\times L_\cx}{\co/U_\cu\times\prod
G_{V^i_m}/U_{\cv^i_m}} & \rTTo^{\pi''_*}_\sim &
\edc{P_\cw\times L_\cu}{\co/U_\cu\times\prod G_{V^i_m}/P_{\cv^i_m}}
& \rTTo_\sim^{\pi'''_*} & \edc{P_\cw}{\co/P_\cx} \\
&& \dTTo<\wr>{\pi''_*} & \ruId & \uTTo<\wr>{\pi'_*}
& \ruTTo^{\pi_*}_\sim \ruTTo(2,4)_{\pi_*}^\sim & \\
\dTTo<{\prod\pi^i_{m*}}>\wr &&
\edc{P_\cw\times L_\cu}{\co/U_\cu\times\prod G_{V^i_m}/P_{\cv^i_m}} &
\rTTo_\sim^{\pi^{\prime*}} &
\edc{P_\cw\times L_\cv\times L_\cu}{\co/U_\cv\times\prod
G_{V^i}/U_{\cv^i}\times\prod G_{V^i_m}/P_{\cv^i_m}}
&& \dTTo>{\beta_{\cx!}} \\
& \ruTTo^{\pr_2^*} && \ruTTo(4,2)_{\pr_3^*} \ruTTo(2,4)^{\pr_{23}^*} &
\dTTo>{\pi_*} && \\
\edc{\prod G_{V^i_m}}{\prod G_{V^i_m}/P_{\cv^i_m}} &&&&
\edc{P_\cw\times L_\cv}{\co/U_\cv\times\prod G_{V^i}/P_{\cy^i}} &&
\edc{P_\cw}{\co/P_\cu} \\
& \rdTTo>{\pr_2^*} &&& \uTTo>{\pr_2^*} & \rdTTo(2,4)^{\beta_!} & \\
\dTTo<{\prod\alpha^i_{m!}} && \edc{\prod G_{V^i}\times L_{\cv^i}}{\prod
G_{V^i}/U_{\cv^i}\times\prod G_{V^i_m}/P_{\cv^i_m}} &
\rTTo^{\pi'''_*}_\sim & \edc{\prod G_{V^i}}{\prod G_{V^i}/P_{\cy^i}} &&
\uTTo>{\pi'''_*}<\sim \\
&& \dTTo>{\pr_{1!}} && \dTTo>{\beta_!} && \\
\edc{\prod G_{V^i_m}}{pt} & \rTTo^{\prod\phi^{i*}} &
\edc{\prod G_{V^i}\times L_{\cv^i}}{\prod G_{V^i}/U_{\cv^i}} &
\rTTo^{\nquad\prod\pi_*^i\nquad}_\sim &
\edc{\prod G_{V^i}}{\prod G_{V^i}/P_{\cv^i}}
& \rTTo^{\pr_2^*\hspace*{1em}} &
\edc{P_\cw\times L_\cv}{\co/U_\cv\times\prod G_{V^i}/P_{\cv^i}}
\end{diagram}
\hspace*{1mm}
}
\caption{Modified first isomorphism}
\label{fig-ll-assoc-2}
\end{figure}

Two lower squares and two rightmost squares in diagrams in
\figref{fig-ll-assoc-2} and \figref{fig-ru-assoc-1} cancel and we have
to prove the following equation.
\begin{figure}
\begin{diagram}[width=4.0em,height=2.400000em]
\edc{\prod G_{V^i_{mp}}}{pt} & \rTTo^{\phi_\cx^*} &
\edc{P_\cw\times L_\cx}{\co/U_\cx} & \rTTo^{\pi_{\cx*}}_\sim &
\edc{P_\cw}{\co/P_\cx} \\
\dTTo<{\prod\phi_m^{i*}} && \uTTo<\wr>{\pi'_*} &
\ruTTo(2,4)^{\pi''_*}_\sim \ruTTo(2,6)_{\pi_*}^\sim & \\
\edc{\prod G_{V^i_m}\times L_{\cv^i_m}}{\prod G_{V^i_m}/U_{\cv^i_m}} &
\rTTo^{\pr_2^*} & \edc{P_\cw\times L_\cu\times
L_\cx}{\co/U_\cu\times\prod G_{V^i_m}/U_{\cv^i_m}} && \\
&& \dTTo<\wr>{\pi''_*} && \uTTo<\wr>{\pi'''_*} \\
\dTTo<{\prod\pi^i_{m*}}>\wr &&
\edc{P_\cw\times L_\cu}{\co/U_\cu\times\prod G_{V^i_m}/P_{\cv^i_m}}
&& \\
& \ruTTo^{\pr_2^*} & \uTTo<\wr>{\pi'_*} && \\
\edc{\prod G_{V^i_m}}{\prod G_{V^i_m}/P_{\cv^i_m}} & \rTTo^{\pr_3^*} &
\edc{P_\cw\times L_\cv\times L_\cu}{\co/U_\cv\times\prod
G_{V^i}/U_{\cv^i}\times\prod G_{V^i_m}/P_{\cv^i_m}} &
\rTTo^{\pi_*}_\sim &
\edc{P_\cw\times L_\cv}{\co/U_\cv\times\prod G_{V^i}/P_{\cy^i}} \\
& \rdTTo^{\pr_2^*} & \uTTo>{\pr_{23}^*} && \uTTo>{\pr_2^*} \\
&& \edc{\prod G_{V^i}\times L_{\cv^i}}{\prod
G_{V^i}/U_{\cv^i}\times\prod G_{V^i_m}/P_{\cv^i_m}} &
\rTTo^{\pi'''_*}_\sim &
\edc{\prod G_{V^i}}{\prod G_{V^i}/P_{\cy^i}}
\end{diagram}
\caption{Part of the first isomorphism}
\label{fig-ll-assoc-3}
\end{figure}
We have to show that the isomorphism in \figref{fig-ll-assoc-3} is
equal to the isomorphism in \figref{fig-ru-assoc-3}.

\begin{figure}
\begin{diagram}[width=4.0em,height=2.400000em]
\edc{\prod G_{V^i_{mp}}}{pt} & \rTTo^{\phi_\cx^*} &
\edc{P_\cw\times L_\cx}{\co/U_\cx} & \rTTo^{\pi_{\cx*}}_\sim &
\edc{P_\cw}{\co/P_\cx} \\
& \rdTTo(2,6)_{\phi^*} \rdTTo(2,4)_{\prod\phi^*_{\cy^i}} &
\uTTo<\wr>{\pi'_*} && \uTTo<\wr>{\pi'''_*} \\
&& \edc{P_\cw\times L_\cv\times L_\cx}{\co/U_\cv\times\prod
G_{V^i}/U_{\cy^i}} & \rTTo^{\pi''_*} &
\edc{P_\cw\times L_\cv}{\co/U_\cv\times\prod G_{V^i}/P_{\cy^i}} \\
\dTTo<{\prod\phi_m^{i*}} && \uTTo>{\pr_2^*} && \\
&& \edc{\prod G_{V^i}\times L_{\cy^i}}{\prod G_{V^i}/U_{\cy^i}} && \\
&& \uTTo<\wr>{\pi'_*} & \rdTTo(2,4)_\sim^{\prod\pi_{\cy^i*}} &
\uTTo>{\pr_2^*} \\
\edc{\prod G_{V^i_m}\times L_{\cv^i_m}}{\prod G_{V^i_m}/U_{\cv^i_m}} &
\rTTo^{\pr_2^*} & \edc{\prod G_{V^i}\times L_{\cv^i}\times
L_{\cy^i}}{\prod G_{V^i}/U_{\cv^i}\times\prod G_{V^i_m}/U_{\cv^i_m}} &&
\\
\dTTo<{\prod\pi^i_{m*}}>\wr && \dTTo<\wr>{\pi''_*} &
\rdTTo_{\pi_*}^\sim & \\
\edc{\prod G_{V^i_m}}{\prod G_{V^i_m}/P_{\cv^i_m}} & \rTTo^{\pr_2^*} &
\edc{\prod G_{V^i}\times L_{\cv^i}}{\prod G_{V^i}/U_{\cv^i}\times\prod
G_{V^i_m}/P_{\cv^i_m}} & \rTTo^{\pi'''_*}_\sim &
\edc{\prod G_{V^i}}{\prod G_{V^i}/P_{\cy^i}}
\end{diagram}
\caption{Part of the second isomorphism}
\label{fig-ru-assoc-3}
\end{figure}

Using \propref{pro-iota-cocycle} we reduce both isomorphisms to the
expressions in \figref{fig-ll-assoc-end} and \figref{fig-ru-assoc-end}.

\begin{figure}
\begin{diagram}[width=4.0em,height=2.400000em]
&& \edc{P_\cw\times L_\cx}{\co/U_\cx} && \\
& \ruTTo^{\phi_\cx^*} & \uTTo<\wr>{\pi_*} & \rdTTo^{\pi_{\cx*}}_\sim &
\\
\edc{\prod G_{V^i_{mp}}}{pt} & \rTTo^{\phi^*} &
\edc{P_\cw\times L_\cv\times L_\cu\times L_\cx}
{\co/U_\cv\times\prod G_{V^i}/U_{\cv^i}\times\prod G_{V^i_m}/U_{\cv^i_m}}
& \rTTo^{\pi_*}_\sim & \edc{P_\cw}{\co/P_\cx} \\
\dTTo<{\prod\phi_m^{i*}} & \ruTTo^{\pr_3^*} & \dTTo<\wr>{\pi_*} &
\ruTTo^{\pi_*}_\sim & \uTTo<\wr>{\pi'''_*} \\
\edc{\prod G_{V^i_m}\times L_{\cv^i_m}}{\prod G_{V^i_m}/U_{\cv^i_m}}
&& \edc{P_\cw\times L_\cv\times L_\cu}{\co/U_\cv\times\prod
G_{V^i}/U_{\cv^i}\times\prod G_{V^i_m}/P_{\cv^i_m}} &
\rTTo^{\pi_*}_\sim &
\edc{P_\cw\times L_\cv}{\co/U_\cv\times\prod G_{V^i}/P_{\cy^i}} \\
\dTTo<{\prod\pi^i_{m*}}>\wr & \ruTTo^{\pr_3^*} & \uTTo>{\pr_{23}^*} &&
\uTTo>{\pr_2^*} \\
\edc{\prod G_{V^i_m}}{\prod G_{V^i_m}/P_{\cv^i_m}} & \rTTo^{\pr_2^*} &
\edc{\prod G_{V^i}\times L_{\cv^i}}{\prod
G_{V^i}/U_{\cv^i}\times\prod G_{V^i_m}/P_{\cv^i_m}} &
\rTTo^{\pi'''_*}_\sim &
\edc{\prod G_{V^i}}{\prod G_{V^i}/P_{\cy^i}}
\end{diagram}
\caption{First simplified isomorphism}
\label{fig-ll-assoc-end}
\end{figure}

\begin{figure}
\begin{diagram}[width=4.0em,height=2.400000em]
&& \edc{P_\cw\times L_\cx}{\co/U_\cx} && \\
& \ruTTo^{\phi_\cx^*} & \uTTo<\wr>{\pi_*} & \rdTTo^{\pi_{\cx*}}_\sim &
\\
\edc{\prod G_{V^i_{mp}}}{pt} & \rTTo^{\phi^*} &
\edc{P_\cw\times L_\cv\times L_\cu\times L_\cx}
{\co/U_\cv\times\prod G_{V^i}/U_{\cv^i}\times\prod G_{V^i_m}/U_{\cv^i_m}}
& \rTTo^{\pi_*}_\sim & \edc{P_\cw}{\co/P_\cx} \\
\dTTo<{\prod\phi_m^{i*}} & \rdTTo^{\phi^*} & \uTTo>{\pr^*_{23}} &
\rdTTo^{\pi_*}_\sim & \uTTo<\wr>{\pi'''_*} \\
\edc{\prod G_{V^i_m}\times L_{\cv^i_m}}{\prod G_{V^i_m}/U_{\cv^i_m}}
& \rTTo^{\pr_2^*} & \edc{\prod G_{V^i}\times L_{\cv^i}\times
L_{\cy^i}}{\prod G_{V^i}/U_{\cv^i}\times\prod G_{V^i_m}/U_{\cv^i_m}} &&
\edc{P_\cw\times L_\cv}{\co/U_\cv\times\prod G_{V^i}/P_{\cy^i}} \\
\dTTo<{\prod\pi^i_{m*}}>\wr && \dTTo<\wr>{\pi''_*} &
\rdTTo^{\pi_*}_\sim & \uTTo>{\pr_2^*} \\
\edc{\prod G_{V^i_m}}{\prod G_{V^i_m}/P_{\cv^i_m}} & \rTTo^{\pr_2^*} &
\edc{\prod G_{V^i}\times L_{\cv^i}}{\prod
G_{V^i}/U_{\cv^i}\times\prod G_{V^i_m}/P_{\cv^i_m}} &
\rTTo^{\pi'''_*}_\sim &
\edc{\prod G_{V^i}}{\prod G_{V^i}/P_{\cy^i}}
\end{diagram}
\caption{Second simplified isomorphism}
\label{fig-ru-assoc-end}
\end{figure}

The two upper triangles in \figref{fig-ll-assoc-end} and
\figref{fig-ru-assoc-end} coincide. Due to \propref{pro-iota-cocycle}
the lower parts are equal to isomorphisms in \figref{fig-assoc-final}.

\begin{figure}
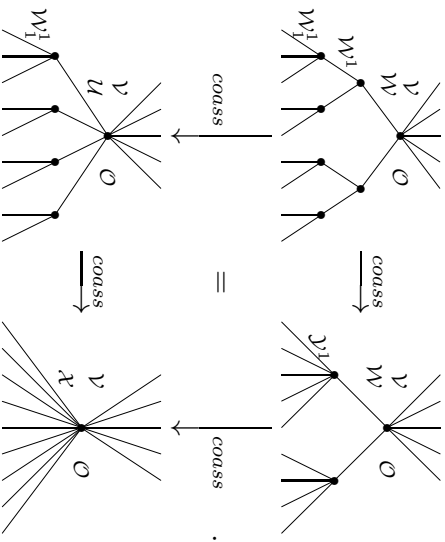

\rotatebox{90}{%
\hspace*{2mm}%
\begin{diagram}[width=4.0em,height=2.600000em,inline]
&&&&&& \\
\edc{\prod G_{V^i_{mp}}}{pt} && \rTTo^{\phi^*} &&
\edc{P_\cw\times L_\cv\times L_\cu\times L_\cx}
{\co/U_\cv\times\prod G_{V^i}/U_{\cv^i}\times\prod G_{V^i_m}/U_{\cv^i_m}}
&& \rTTo^{\pi_*}_\sim && \edc{P_\cw}{\co/P_\cx} \\
\dTTo<{\prod\phi_m^{i*}} &&& \ruTTo(4,2)^{\pr_3^*} \ruTTo^{\pr_{23}^*} &&
\rdTTo<\wr>{\pi_*} && \ruTTo^{\pi'_*}_\sim & \uTTo<\wr>{\pi'''_*} \\
\edc{\prod G_{V^i_m}\times L_{\cv^i_m}}{\prod G_{V^i_m}/U_{\cv^i_m}} &
\rTTo^{\pr_2^*} & \edc{\prod G_{V^i}\times L_{\cv^i}\times
L_{\cy^i}}{\prod G_{V^i}/U_{\cv^i}\times\prod G_{V^i_m}/U_{\cv^i_m}}
\hspace*{-2cm} &&&& \hspace*{-2cm} \edc{P_\cw\times L_\cv\times L_\cu}
{\co/U_\cv\times\prod G_{V^i}/U_{\cv^i}\times\prod G_{V^i_m}/P_{\cv^i_m}}
& \rTTo^{\pi_*}_\sim &
\edc{P_\cw\times L_\cv}{\co/U_\cv\times\prod G_{V^i}/P_{\cy^i}} \\
\dTTo<{\prod\pi^i_{m*}}>\wr &&& \rdTTo^{\pi''_*} && \ruTTo>{\pr_{23}^*}
&&& \uTTo>{\pr_2^*} \\
\edc{\prod G_{V^i_m}}{\prod G_{V^i_m}/P_{\cv^i_m}} && \rTTo^{\pr_2^*}
&& \edc{\prod G_{V^i}\times L_{\cv^i}}{\prod
G_{V^i}/U_{\cv^i}\times\prod G_{V^i_m}/P_{\cv^i_m}} &&
\rTTo^{\pi'''_*}_\sim && \edc{\prod G_{V^i}}{\prod G_{V^i}/P_{\cy^i}}
%%%
\hspace*{1mm} \\
\hspace*{1mm} \\
%%% next diagram
\edc{\prod G_{V^i_{mp}}}{pt} && \rTTo^{\phi^*} &&
\edc{P_\cw\times L_\cv\times L_\cu\times L_\cx}
{\co/U_\cv\times\prod G_{V^i}/U_{\cv^i}\times\prod G_{V^i_m}/U_{\cv^i_m}}
&& \rTTo^{\pi_*}_\sim && \edc{P_\cw}{\co/P_\cx} \\
\dTTo<{\prod\phi_m^{i*}} & \rdTTo^{\pr_3^*} && \ruTTo^{\pr_{23}^*} &&
\rdTTo<\wr>{\pi_*} \rdTTo(4,2)^{\pi_*}_\sim &&& \uTTo<\wr>{\pi'''_*} \\
\edc{\prod G_{V^i_m}\times L_{\cv^i_m}}{\prod G_{V^i_m}/U_{\cv^i_m}} &
\rTTo^{\pr_2^*} & \edc{\prod G_{V^i}\times L_{\cv^i}\times
L_{\cy^i}}{\prod G_{V^i}/U_{\cv^i}\times\prod G_{V^i_m}/U_{\cv^i_m}}
\hspace*{-2cm} &&&& \hspace*{-2cm} \edc{P_\cw\times L_\cv\times L_\cu}
{\co/U_\cv\times\prod G_{V^i}/U_{\cv^i}\times\prod G_{V^i_m}/P_{\cv^i_m}}
& \rTTo^{\pi_*}_\sim &
\edc{P_\cw\times L_\cv}{\co/U_\cv\times\prod G_{V^i}/P_{\cy^i}} \\
\dTTo<{\prod\pi^i_{m*}}>\wr &&& \rdTTo^{\pi''_*} && \ruTTo>{\pr_{23}^*}
&&& \uTTo>{\pr_2^*} \\
\edc{\prod G_{V^i_m}}{\prod G_{V^i_m}/P_{\cv^i_m}} && \rTTo^{\pr_2^*}
&& \edc{\prod G_{V^i}\times L_{\cv^i}}{\prod
G_{V^i}/U_{\cv^i}\times\prod G_{V^i_m}/P_{\cv^i_m}} &&
\rTTo^{\pi'''_*}_\sim && \edc{\prod G_{V^i}}{\prod G_{V^i}/P_{\cy^i}}
\end{diagram}
\hspace*{5mm}
}
\caption{Final isomorphisms}
\label{fig-assoc-final}
\end{figure}

Finally, the last two isomorphisms are equal to each other due again to
\propref{pro-iota-cocycle}.
\end{proof}

\subsection{Comultiplication}
The comultiplication functor is
\begin{equation*}
\qquad
\vstretch 60
\begin{tangle}
\step[4]\idash\noderu{\cw} \\
\nodeld{W^1}\sw4\sw0\nodeld{W^2}\sw2\sw0\s\se0\se2\se0\se4\noderd{W^l}
\end{tangle}
\quad = \pitchfork_\cw^\cv =
\Bigl( \edc{P_\cw}{pt} \rTTo^{\Resto_{L_\cw,P_\cw}} \edc{L_\cw}{pt} \Bigr)
= i_{L_\cw,P_\cw}^*.
\end{equation*}
Together with multiplication definition of \secref{sec-Multiplication}
it gives general operation (in the set-up of
\eqref{eq-decomp-V-to-Vi2}--\eqref{eq-decomp-V-to-Wj2})
\begin{multline*}
\zhe_\cw^{\co;\cv} = \qquad
\vstretch 80
\begin{tangle}
\nodelu{V^1}\Step\nodelu{V^2}\step[6]\noderu{V^k} \\
\nw4\nw0\nw2\nodel{\co}\nw0\n\ne0\ne2\ne0\ne4 \\
\nodeld{W^1}\sw4\sw0\nodeld{W^2}\sw2\sw0\id\se0\se2\se0\se4\noderd{W^l}
\end{tangle}
\qquad = \Bigl( D_{L_\cv}(pt) \rTTo^{\phi^*}
D_{P_\cw\times L_\cv}(\co/{U_\cv}) \\[2mm]
\rTTo^{\pi_*} D_{P_\cw}(\co/{P_\cv}) \rTTo^{\alpha_!} D_{P_\cw}(pt)
\rTTo^{\Resto_{L_\cw,P_\cw}} D_{L_\cw}(pt) \Bigr)
\end{multline*}
\[
\simeq \bigl(  D_{L_\cv}(pt) \rTTo^{\phi^*}
D_{L_\cw\times L_\cv}(\co/{U_\cv}) \rTTo^{\pi_*} D_{L_\cw}(\co/{P_\cv})
\rTTo^{\alpha_!} D_{L_\cw}(pt) \bigr).
\]

In the particular case $k=1$ we have $\co=L_\cv=P_\cv=G_V$, $U_\cv=1$,
and the comultiplication operation is isomorphic to
\begin{multline*}
\vstretch 80
\begin{tanglec}
\n \\
\noded{W^1}\sw4\sw0\noded{W^2}\sw2\sw0\id\se0\se2\se0\se4\noded{W^l}
\end{tanglec}
\simeq \bigl( D_{G_V}(pt) \rTTo^{\Resto_{P_\cw,G_V}} D_{P_\cw}(pt)
\rTTo^{\Resto_{L_\cw,P_\cw}} D_{L_\cw}(pt) \bigr) \\
\simeq \bigl( D_{G_V}(pt) \rTTo^{\Resto_{L_\cw,G_V}} D_{L_\cw}(pt) \bigr)
\end{multline*}
by \propref{pro-psi1=res}.

\subsection{The coassociativity isomorphism}
Assume that besides decomposition \eqref{eq-decomp-V-to-Wj2} of $V$
into direct sum of $W^j$ we have also decompositions
\[ \cw^j:\qquad W^j = \oplus_{m=1}^{m^j} W_m^j .\]
We can produce out of these the decomposition
\[ \cu:\qquad V = \oplus_{j,m} W_m^j ,\]
where the lexicographic order of summands is used
$W_1^1$, \dots,$W_{m^1}^1,W_1^2$, \dots,$W_{m^2}^2$,
\dots,$W_1^l$, \dots,$W_{m^l}^l$.

Coassociativity isomorphism
$coass:\pitchfork_\cu^\cw\circ\zhe_\cw^\cv\to\zhe_\cu^\cv$ between
functors
\[
\quad
\hstretch 90
\vstretch 80
\begin{tanglec}
\nw2\nw1\id\ne1\noderu{\cv}\ne2 \\
\nodeld{\cw^1\;}\sw3\Step\nodeld{\cw^2\;}\s\Step\se3\noderd{\cw^l} \\
\sw1\s\se1\step\sw1\s\se1\step\sw1\s\se1
\end{tanglec}
\;= \Bigl( \edc{L_\cv}{pt} \rTTo^{\Psi^\cv_\cw} \edc{P_\cw}{pt}
\rTTo^{i^*_{L_\cw,P_\cw}} \edc{L_\cw}{pt} \rTTo^{i^*_{L_\cu,L_\cw}}
\edc{L_\cu}{pt} \Bigr)
\]
\[
\text{and}\quad
\vstretch 100
\begin{tanglec}
\nw2\nw1\id\ne1\noderd{\CU}\ne2 \\
\sw4\sw3\sw2\sw1\s\se1\se2\se3\se4
\end{tanglec}
\quad = \Bigl( \Dbc_{L_\cv}(pt) \rTTo^{\Psi^\cv_\cu} \Dbc_{P_\cu}(pt)
\rTTo^{i^*_{L_\cu,P_\cu}} \Dbc_{L_\cu}(pt) \Bigr)
\]
is given by the diagram
\begin{diagram}[height=2.7em,nobalance]
\edc{L_\cv}{pt} & \rTTo^{\phi^*} & \edc{P_\cw\times L_\cv}{\co/{U_\cv}} &
\rTTo^{\pi_*} & \edc{P_\cw}{\co/{P_\cv}} & \rTTo^{\alpha_!} &
\edc{P_\cw}{pt} & \rTTo^{i^*_{L_\cw,P_\cw}} & \edc{L_\cw}{pt} \\
\dId && \dTTo~{i^*_{P_\cu\times L_\cv,P_\cw\times L_\cv}} &&
\dTTo~{i^*_{P_\cu,P_\cw}} && \dTTo~{i^*_{P_\cu,P_\cw}} &
\rdTTo^{i^*_{L_\cu,P_\cw}} & \dTTo~{i^*_{L_\cu,L_\cw}} \\
\edc{L_\cv}{pt} & \rTTo^{\phi^*} & \edc{P_\cu\times L_\cv}{\co/{U_\cv}} &
\rTTo^{\pi_*} & \edc{P_\cu}{\co/{P_\cv}} & \rTTo^{\alpha_!} &
\edc{P_\cu}{pt} & \rTTo^{i^*_{L_\cu,P_\cu}} & \edc{L_\cu}{pt}
\end{diagram}

\subsection{The coassociativity equation}
Assume that we have the following decompositions
of finite-dimensional $\CC$\n-vector spaces:
\begin{alignat*}4
\cv &:&\quad V &= \oplus_i V^i, &\qquad
\cw &:&\quad V &= \oplus_j W^j, \\
\cw^j &:&\quad W^j &= \oplus_m W^j_m, &\qquad
\cw^j_m &:&\quad W^j_m &= \oplus_p W^j_{m,p}.
\end{alignat*}
These decompositions imply also the following decompositions:
\begin{equation*}
\cu :\; V = \oplus_{j,m} W^j_m, \quad\;
\cy^j :\; W^j = \oplus_{m,p} W^j_{m,p}, \quad\;
\cx :\; V = \oplus_{j,m,p} W^j_{m,p}.
\end{equation*}
Let $\co\subset G_V$ be a left $P_\cw$\n-invariant
and right $P_\cv$\n-invariant subset.
Coassociativity isomorphisms give two isomorphisms
between the composite operation and the single operation
as in the diagram
\begin{diagram}[height=4.5em,LaTeXeqno]
\vstretch 75
\begin{tanglec}
\nw2\nodelu{\cv}\nodeld{\cw}\nw1\n\ne1\noder{\co}\ne2 \\
\nodeld{\cw^1}\node\sw2\Step\se2\node \\
\nodeld{\cw^1_1}\sw1\se1\Step\sw1\se1 \\
\sw1\s\se1\step\s\se1\step\s\se1\step\s\se1
\end{tanglec}
& \rTTo^{coass} &
\begin{tanglec}
\nw2\nodelu{\cv}\nodeld{\cw}\nw1\n\ne1\noder{\co}\ne2 \\
\nodeld{\cy^1}\sw2\Step\se2 \\
\sw2\sw1\s\se1\se2\step\sw1\s\se1\se2
\end{tanglec}
\\
\dTTo<{coass} & = & \dTTo>{coass} & . \\
\begin{tanglec}
\nw2\nodelu{\cv}\nodeld{\CU}\nw1\n\ne1\noder{\co}\ne2 \\
\nodeld{\cw^1_1}\sw3\sw0\sw1\se1\se0\se3 \\
\sw1\s\se1\step\s\se1\step\s\se1\step\s\se1
\end{tanglec}
& \rTTo^{coass} &
\vstretch 150
\begin{tanglec}
\nw2\nodelu{\cv}\nodeld{\cx}\nw1\n\ne1\noder{\co}\ne2 \\
\sw4\sw3\sw2\sw1\id\se1\se2\se3\se4
\end{tanglec}
\label{dia-coass-square}
\end{diagram}

\begin{proposition}[The coassociativity equation]
Diagram \eqref{dia-coass-square} is commutative.
\end{proposition}

\begin{proof}
The graphically expressed equation means equality
of the following two isomorphisms
\begin{diagram}[nobalance,width=2.5em]
\edc{L_\cv}{pt} & \rTTo^{\phi^*} & \edc{P_\cw\times L_\cv}{\co/{U_\cv}} &
\rTTo^{\pi_*} & \edc{P_\cw}{\co/{P_\cv}} & \rTTo^{\alpha_!} &
\edc{P_\cw}{pt} & \rTTo^{i^*_{L_\cw,P_\cw}} & \edc{L_\cw}{pt} \\
\dId && \dTTo>{i^*_{P_\cu\times L_\cv,P_\cw\times L_\cv}} &&
\dTTo>{i^*_{P_\cu,P_\cw}} && \dTTo>{i^*_{P_\cu,P_\cw}} &
\rdTTo^{i^*_{L_\cu,P_\cw}} & \dTTo>{i^*_{L_\cu,L_\cw}} \\
\edc{L_\cv}{pt} & \rTTo^{\phi^*} & \edc{P_\cu\times L_\cv}{\co/{U_\cv}} &
\rTTo^{\pi_*} & \edc{P_\cu}{\co/{P_\cv}} & \rTTo^{\alpha_!} &
\edc{P_\cu}{pt} & \rTTo^{i^*_{L_\cu,P_\cu}} & \edc{L_\cu}{pt} \\
\dId && \dTTo>{i^*_{P_\cx\times L_\cv,P_\cu\times L_\cv}} &&
\dTTo>{i^*_{P_\cx,P_\cu}} && \dTTo>{i^*_{P_\cx,P_\cu}} &
\rdTTo^{i^*_{L_\cx,P_\cu}} & \dTTo>{i^*_{L_\cx,L_\cu}} \\
\edc{L_\cv}{pt} & \rTTo^{\phi^*} & \edc{P_\cx\times L_\cv}{\co/{U_\cv}} &
\rTTo^{\pi_*} & \edc{P_\cx}{\co/{P_\cv}} & \rTTo^{\alpha_!} &
\edc{P_\cx}{pt} & \rTTo^{i^*_{L_\cx,P_\cx}} & \edc{L_\cx}{pt}
\end{diagram}
\begin{diagram}[nobalance,width=3.0em,tight]
\edc{L_\cv}{pt} & \rTTo^{\phi^*\hspace*{-1mm}} &
\edc{P_\cw\times L_\cv}{\co/{U_\cv}} & \rTTo^{\hspace*{-1mm}\pi_*} &
\edc{P_\cw}{\co/{P_\cv}} & \rTTo^{\alpha_!} &
\edc{P_\cw}{pt} & \rTTo^{i^*_{L_\cw,P_\cw}} & \edc{L_\cw}{pt} \\
&&& &&& &&& \rdTTo>{i^*_{L_\cu,L_\cw}} & \\
\dId &&
\dTTo~{\hspace*{-1mm}i^*_{P_\cx\times L_\cv,P_\cw\times L_\cv}\hspace*{-1mm}}
&& \dTTo~{i^*_{P_\cx,P_\cw}} && \dTTo~{i^*_{P_\cx,P_\cw}} &&
\dTTo~{i^*_{L_\cx,L_\cw}} && \edc{L_\cu}{pt}  \\
&&& &&& &&& \ldTTo>{i^*_{L_\cx,L_\cu}} & \\
\edc{L_\cv}{pt} & \rTTo^{\phi^*\hspace*{-1mm}} &
\edc{P_\cx\times L_\cv}{\co/{U_\cv}} & \rTTo^{\hspace*{-1mm}\pi_*} &
\edc{P_\cx}{\co/{P_\cv}} & \rTTo^{\alpha_!} &
\edc{P_\cx}{pt} & \rTTo^{i^*_{L_\cx,P_\cx}} & \edc{L_\cx}{pt}
\end{diagram}
Four triangular prisms, which follow from \propref{pro-iota-cocycle}
and \propref{prop-distr-base-chge},
imply that the above isomorphisms are equal.
\end{proof}

\subsection{The coherence isomorphism}
Assume given a vector space decomposition
\[ V = \oplus_{m=1}^k \oplus_{r=1}^l V_r^m .\]
We shall denote it $\cy$ if the summands are ordered as follows
\[ V_1^1,V_2^1,\dots,V^1_l,V^2_1,V^2_2,\dots,V^2_l,\dots,
V_1^k,V_2^k,\dots,V_l^k .\]
The same decomposition with order of summands
\[ V_1^1,V_1^2,\dots,V_1^k,V_2^1,V_2^2,\dots,V_2^k,\dots,
V_l^1,V_l^2,\dots,V_l^k \]
will be denoted $\cx$. We also have decompositions
\begin{alignat*}4
\cv^m &:&\quad  V^m &= \oplus_{r=1}^l V^m_r &\qquad
\cv &:&\quad  V &= \oplus_{m=1}^k V^m \\
\cw_r &:&\quad  W_r &= \oplus_{m=1}^k V_r^m &\qquad
\cw &:&\quad  V &= \oplus_{r=1}^l W_r
\end{alignat*}
which produce filtrations $V^m_{(r)}$ of $V^m$, $V^{(m)}$ of $V$,
$V_r^{(m)}$ of $W_r$ and $W_{(r)}$ of $V$. We associate with these
filtrations parabolic groups $P_{\cv^m}\subset G_{V^m}$,
$P_\cv\subset G_V$, $P_{\cw_r}\subset G_{W_r}$ and $P_\cw\subset G_\cw$.

\begin{notations}
Let $A$, $B$ be vector spaces. We denote the braiding functor as
\begin{equation*}
%\vstretch 130
\begin{tanglec}
\xx
\end{tanglec}
= \bigl( D_{G_{A}\times G_{B}}(pt) \xra\tau D_{G_{A}\times G_{B}}(pt)
\xra{\sigma^*} D_{G_{B}\times G_{A}}(pt) \bigr),
\end{equation*}
where $\sigma$ is the permutation isomorphism of groups and modules and
the functor $\tau$ is the shift \( \tau(L) = L[-2\dim A\dim B] \).
\end{notations}

We claim that for any collection of indices and
for any collection of bi-invariant locally closed subsets
\((\co'_1,\dots,\co'_k,\co''_1,\dots,\co''_l)\),
which may occur in the following diagram,
there exists a bi-invariant locally closed subset $\co$
and a coherence isomorphism
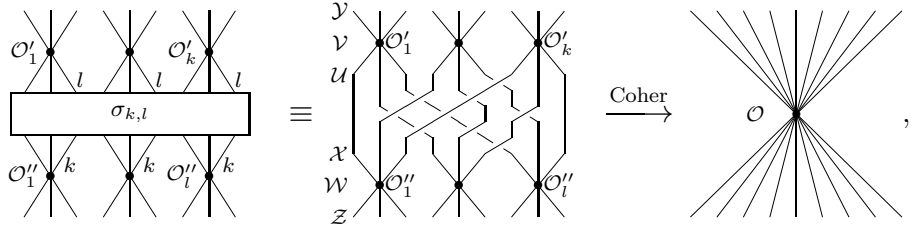
\[
\vstretch 78
\begin{tanglec}
\nw1\nodel{\co'_1\;}\n\ne1\step\nw1\n\ne1\step\nw1\nodel{\co'_k\;}\n\ne1 \\
\sw1\s\se1\nodeu{\;l}\step\sw1\s\se1\nodeu{\;l}\step\sw1\s\se1\nodeu{\;l} \\
\ffbox{9}{\sigma_{k,l}} \\
\nw1\nodel{\co''_1\;}\n\noderu{k}\ne1\step\nw1\n\noderu{k}\ne1\step%
\nw1\nodel{\co''_l\;}\n\noderu{k}\ne1 \\
\sw1\s\se1\step\sw1\s\se1\step\sw1\s\se1
\end{tanglec}
\quad\equiv\quad
\vstretch 60
\begin{tangle}
\nodel{\cy} \\
\nodel{\cv}\nw1\noder{\co'_1}\n\ne1\step\nw1\n\ne1\step\nw1\noder{\co'_k}\n\ne1 \\
\nodel{\CU}\sw1\s\se1\step\sw1\s\se1\step\sw1\s\se1 \\
\hh \id\step\id\step\id\step\id\step\id\step\xx\step\id\step\id \\
\hh \id\step\id\step\xx\step\xx\step\id\step\id\step\id \\
\hh \id\step\xx\step\xx\step\id\step\xx\step\id \\
\hh \id\step\id\step\xx\step\xx\step\id\step\id\step\id \\
\hh \nodel{\cx}\id\step\id\step\id\step\id\step\id\step\xx\step\id\step\id \\
\nodel{\cw}\nw1\noder{\co''_1}\n\ne1\step\nw1\n\ne1\step\nw1\noder{\co''_l}\n\ne1 \\
\nodel{\cz}\sw1\s\se1\step\sw1\s\se1\step\sw1\s\se1
\end{tangle}
\quad \rTTo^{\Coher} \;
\vstretch 195
\begin{tangle}
\nw4\nw3\nw2\nodel{\co}\nw1\n\ne1\ne2\ne3\ne4 \\
\sw4\sw3\sw2\sw1\id\se1\se2\se3\se4
\end{tangle}
,\]
explicitly described below.
Here $\sigma_{k,l}=(s_{k,l})_+\sptilde$ is the braid, corresponding
to the permutation $s_{k,l}$ of the set $\{1,2,\dots,kl\}$,
\begin{equation}
s_{k,l}(1+t+nl)=1+n+tk \text{ \ for \ } 0\le t<l, 0\le n<k,
\label{eq-s-permute}
\end{equation}
under the standard splitting $S_{kl}\to B_{kl}$, which maps
the elementary transpositions to the generators of the braid group.
The subset $\co$ is computed as follows.

Bi-invariance of initial parametrising subsets means that
\begin{align*}
G_{V^m} \supset \co'_m &\in P_{V_-^m}\sets P_{\cy^m}, \\
G_{W_r} \supset \co''_r &\in P_{\cz_r}\sets P_{V_r^-},
\end{align*}
where $\cy^m$ is the incoming filtration of $V^m$, and $\cz_r$ is the
outgoing filtration of $W_r$. The subsets
\begin{align*}
P_{\cv}\supset \overline{\co} &\overset{\text{def}}=
U_\cv \cdot \prod_m\co'_m = \prod_m\co'_m \cdot U_\cv
\in P_{V_-^=}\sets P_{\cy}, \\
P_{\cv}\supset \underline{\co} &\overset{\text{def}}=
U_\cw \cdot \prod_r\co''_r = \prod_r\co''_r \cdot U_\cw
\in P_{\cz}\sets P_{V_=^-}
\end{align*}
are locally closed. The subset
\begin{gather*}
G_V\supset \co \overset{\text{def}}= \und\co\cdot\overline\co
\in P_{\cz}\sets P_{\cy}, \quad
\co = \und\co P_\cu\times_{P_\cu}\overline\co =
\und\co\times_{P_\cw\cap P_\cv}\overline\co,
\end{gather*}
is also locally closed. Indeed, $\und\co\times\overline\co$ is locally
closed in $P_\cw\times P_\cv$, hence,
$\und\co\times_{P_\cw\cap P_\cv}\overline\co$ is locally closed in
$P_\cw\times_{P_\cw\cap P_\cv}P_\cv\simeq P_\cw\cdot P_\cv$.
The subspace $C=P_\cw\cdot P_\cv$ of $G_V$ is locally closed, since it
is an orbit of $P_\cw\times P_\cv$ in $G_V$, and the number of such
orbits is finite. Indeed, $C$ is embedded in
its closure as
\[ C = \overline{C} - \bigcup\{\overline{\co_a} \mid
\co_a\subset\overline{C} \text{ is an orbit of } P_\cw\times P_\cv,
\; \co_a\neq C \} \]
due to dimension considerations.

Both associativity isomorphism and coassociativity isomorphism
are particular cases of the general coherence isomorphism.

In the particular case, when the upper row of operations consists of
comultiplications (operations with single input), and the lower row of
operations consists of multiplications (operations with single output),
we have necessarily $\co'_m=G_{V^m}$ and $\co''_r=G_{W_r}$. In such
cases we omit the parametrising set, since it is unique. Graphical
notation here is the following:
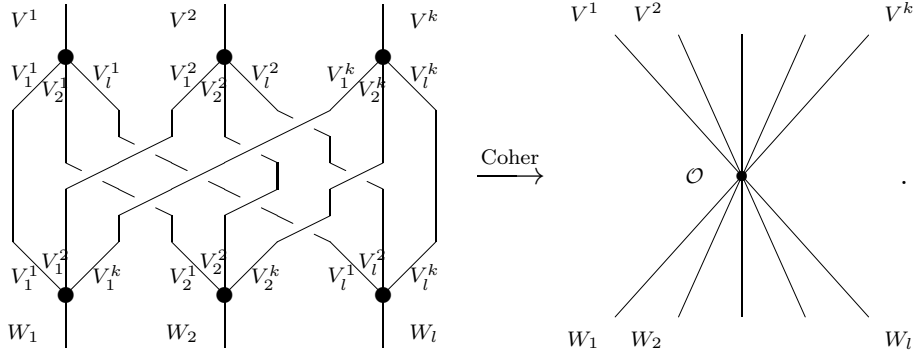
\[
\vstretch 100
\hstretch 200
\begin{tanglec}
\nodeld{V^1}\step[3]\nodeld{V^2}\step[3]\noderd{V^k} \\
\nodeld{V_1^1}\id\noderd{V_l^1}\step[3]\nodeld{V_1^2}\id%
\noderd{V_l^2}\step[3]\nodeld{V_1^k}\id\noderd{V_l^k} \\
\sw1\nodeu{V_2^1\ \ }\s\se1\step\sw1%
\nodeu{V^2_2\ \ }\s\se1\step\sw1\nodeu{V_2^k\ \ }\s\se1 \\
\hh \id\step\id\step\id\step\id\step\id\step\xx\step\id\step\id \\
\hh \id\step\id\step\xx\step\xx\step\id\step\id\step\id \\
\hh \id\step\xx\step\xx\step\id\step\xx\step\id \\
\hh \id\step\id\step\xx\step\xx\step\id\step\id\step\id \\
\hh \id\step\noded{V_1^2\ \ }\id\step\id\step\id\step%
\noded{V_2^2\ \ }\id\step\xx\step\noded{V_l^2\ \ }\id\step\id \\
\nw1\nodelu{V_1^1}\n\noderu{V_1^k}\ne1\step\nw1\nodelu{V_2^1}\n%
\noderu{V_2^k}\ne1\step\nw1\nodelu{V_l^1}\n\noderu{V_l^k}\ne1 \\
\nodelu{W_1}\id\step[3]\nodelu{W_2}\id\step[3]\id\noderu{W_l}
\end{tanglec}
\quad \rTTo^\Coher \qquad
\vstretch 267
\hstretch 120
\begin{tangle}
\nodelu{V^1}\Step\nodelu{V^2}\step[6]\noderu{V^k} \\
\nw4\nw0\nw2\nodel{\co}\nw0\n\ne0\ne2\ne0\ne4 \\
\nodeld{W_1}\sw4\sw0\nodeld{W_2}\sw2\sw0\id\se0\se2\se0\se4\noderd{W_l}
\end{tangle}
\quad. \]
By general formulas we find here $\overline{\co}=P_\cv$,
$\underline{\co}=P_\cw$, hence, $\co=P_\cv P_\cw$.

The general coherence isomorphism is built as the composition
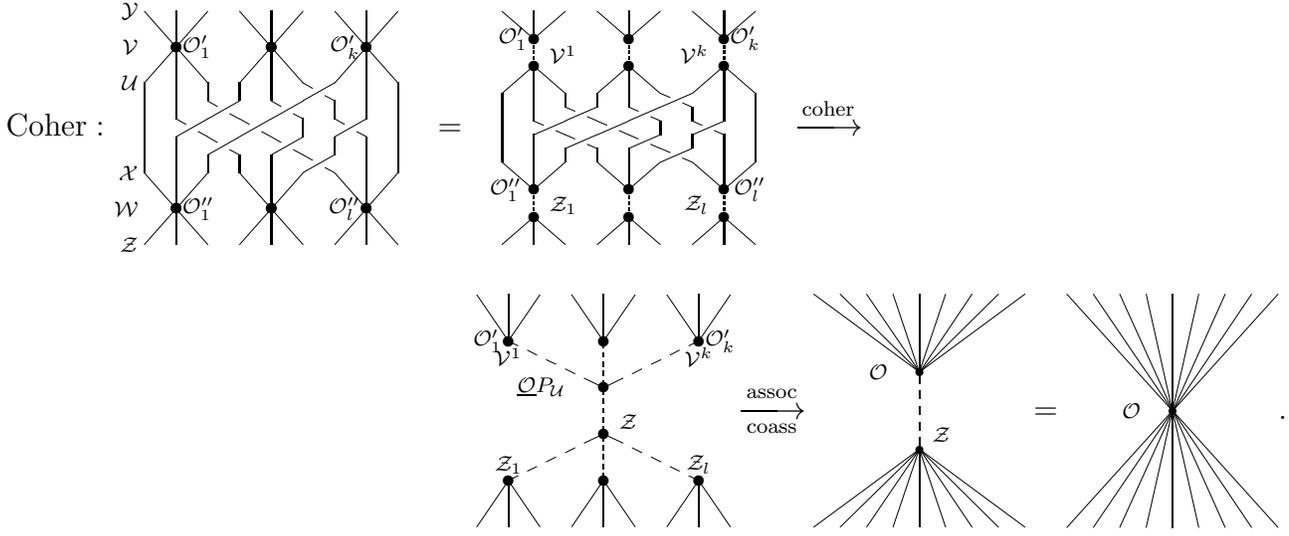
\begin{multline*}
\Coher:\quad
\hstretch 120
\vstretch 68
\begin{tangle}
\nodel{\cy} \\
\nodel{\cv}\nw1\noder{\co'_1}\n\ne1\step\nw1\n\ne1\step\nw1\nodel{\co'_k}\n\ne1 \\
\nodel{\CU}\sw1\s\se1\step\sw1\s\se1\step\sw1\s\se1 \\
\hh \id\step\id\step\id\step\id\step\id\step\xx\step\id\step\id \\
\hh \id\step\id\step\xx\step\xx\step\id\step\id\step\id \\
\hh \id\step\xx\step\xx\step\id\step\xx\step\id \\
\hh \id\step\id\step\xx\step\xx\step\id\step\id\step\id \\
\hh \nodel{\cx}\id\step\id\step\id\step\id\step\id\step\xx\step\id\step\id \\
\nodel{\cw}\nw1\noder{\co''_1}\n\ne1\step\nw1\n\ne1\step\nw1\nodel{\co''_l}\n\ne1 \\
\nodel{\cz}\sw1\s\se1\step\sw1\s\se1\step\sw1\s\se1
\end{tangle}
\quad=\quad
\vstretch 52
\begin{tangle}
\nw1\nodel{\co'_1}\n\ne1\step\nw1\n\ne1\step\nw1\noder{\co'_k}\n\ne1 \\
\step\idash\noderu{\cv^1}\step[3]\idash\step[3]\idash\nodelu{\cv^k} \\
\sw1\s\se1\step\sw1\s\se1\step\sw1\s\se1 \\
\hh \id\step\id\step\id\step\id\step\id\step\xx\step\id\step\id \\
\hh \id\step\id\step\xx\step\xx\step\id\step\id\step\id \\
\hh \id\step\xx\step\xx\step\id\step\xx\step\id \\
\hh \id\step\id\step\xx\step\xx\step\id\step\id\step\id \\
\hh \id\step\id\step\id\step\id\step\id\step\xx\step\id\step\id \\
\nw1\nodel{\co''_1\;}\n\ne1\step\nw1\n\ne1\step\nw1\noder{\,\co''_l}\n\ne1 \\
\step\idash\noderu{\cz_1}\step[3]\idash\step[3]\idash\nodelu{\cz_l} \\
\sw1\s\se1\step\sw1\s\se1\step\sw1\s\se1
\end{tangle}
\quad \rTTo^{\coher} \\[5mm]
\hstretch 120
\vstretch 88
\begin{tangle}
\nw1\nodel{\co'_1}\n\noded{\cv^1}\ne1\step\nw1\n\ne1%
\step\nw1\noder{\co'_k}\n\noded{\cv^k}\ne1 \\
\step\nwd3\step\nodel{\und\co P_\CU}\step\nd\Step\ned3 \\
\step[4]\idash\noderu{\cz} \\
\step\nodeu{\cz_1}\swd3\Step\sd\Step\sed3\nodeu{\cz_l} \\
\sw1\s\se1\step\sw1\s\se1\step\sw1\s\se1
\end{tangle}
\rTTo^{\assoc}_{\coass}
\hstretch 100
\vstretch 147
\begin{tangle}
\nw4\nw3\nw2\nodel{\co}\nw1\n\ne1\ne2\ne3\ne4 \\
\step[4]\idash\noderu{\cz} \\
\sw4\sw3\sw2\sw1\s\se1\se2\se3\se4
\end{tangle}
=
\vstretch 221
\begin{tangle}
\nw4\nw3\nw2\nodel{\co}\nw1\n\ne1\ne2\ne3\ne4 \\
\sw4\sw3\sw2\sw1\id\se1\se2\se3\se4
\end{tangle}
    .
\end{multline*}

The three components of the coherence isomorphism are defined next.
\begin{diagram}[height=2.6em,nobalance]
\Dec{pt}{\prod G_{Y_s^m}} & \rTTo^{\phi^*} &
\Dec{\prod \co'_m}{\prod P_{\cv^m}\times P_{\cy^m}}
& \rTTo^{\pi_*} & \Dec{\prod \co'_m/{P_{\cy^m}}}{\prod P_{\cv^m}}
& \rTTo^{\alpha_!} & \Dec{pt}{\prod P_{\cv^m}} \\
\dTTo<{\phi^*} && \dTTo>{\phi^*} && \dTTo>{\phi^*} && \dTTo>{\phi^*} \\
\Dec{\co}{P_\cz\times P_\cy} & \rTTo^{\pi^*} &
\Dec{\und\co P_\cu\times\overline{\co}}{P_\cz\times P_\cu\times P_\cy}
& \rTTo^{\pi_*} &
\Dec{\und\co P_\cu\times\overline{\co}/{P_\cy}}{P_\cz\times P_\cu}
& \rTTo^{1\times\alpha_!} & \Dec{\und\co P_\cu}{P_\cz\times P_\cu} \\
& \rdId<{\Id} & \dTTo>{\pi_*} && \dTTo>{\pi_*} && \dTTo>{\pi_*} \\
&& \Dec{\co}{P_\cz\times P_\cy} & \rTTo^{\pi_*}
& \Dec{\co/{P_\cy}}{P_\cz} & \rTTo^{\beta_!}
& \Dec{\und\co P_\cu/{P_\cu}}{P_\cz} \\
&& \assoc &&& \rdTTo_{\alpha_!} & \dTTo>{\alpha_!} \\
&&&&&& \Dec{pt}{P_\cz}
\end{diagram}
\[ \coass: \Bigl(\Dec{pt}{P_\cz} \rTTo^{\iota^*}
\Dec{pt}{\prod_rP_{\cz_r}} \rTTo^{\iota^*}
\Dec{pt}{\prod_{n,r}G_{Z^n_r}}\Bigr)
\rTTo^\sim \Bigl(\Dec{pt}{P_\cz} \rTTo^{\iota^*}
\Dec{pt}{\prod_{n,r}G_{Z^n_r}}\Bigr). \]

\begin{figure}
\resizebox{\textwids}{!}{%Sl7coh-4
\begin{diagram}[height=2.5em,width=2.5em,inline,nobalance]
\dec{pt}{\prod_mP_{\cv^m}} &&& \rTTo^{\iota^*}
&&& \dec{pt}{\prod_{m,r}G_{V^m_r}} \\
& \rdId^{\Id} & (1) &&& \ldTTo(4,2)^{p^*} \ldTTo>{p^*}
& \dTTo>{\phi^*} \\
\dTTo<{\phi^*} && \dec{pt}{\prod_mP_{\cv^m}}
&& \dec{\und\co/{U_\cw\cap P_\cv}}{\prod_rP_{\cz_r}\times P_{\cw_r}}
&& \dec{\prod_r\co_r''}{\prod_rP_{\cz_r}\times P_{\cw_r}} \\
&& \dTTo^{p^*} & \ldTTo^{q^*} & \dId<{\ttt(2)\hspace*{0.7em}}>{\Id}
& \ldTTo^{\rho^*} & \dTTo<{\ttt(3)\hspace*{0.5em}}>{T^{-2A}} \\
\dec{\und\co P_\cu}{P_\cz\times P_\cu} & \rTTo^{i^*} &
\dec{\und\co}{(\prod P_{\cz_r})\times(P_\cw\cap P_\cv)}
& \rTTo^{q_*} &
\dec{\und\co/{U_\cw\cap P_\cv}}{\prod_rP_{\cz_r}\times P_{\cw_r}}
& \rTTo^{\rho_!} &
\dec{\prod_r\co_r''}{\prod_rP_{\cz_r}\times P_{\cw_r}} \\
\dTTo<{\pi_*} & (4) & \dTTo<{q_*} & \ldTTo>{q_*} & (5)
&& \dTTo>{\pi_*} \\
\dec{\und\co P_\cu/{P_\cu}}{P_\cz} & \rTTo^{\iota^*}
& \dec{\und\co/{P_\cw\cap P_\cv}}{\prod P_{\cz_r}} && \rTTo^{\rho_!}
&& \dec{\prod_r\co_r''/{P_{\cw_r}}}{\prod_rP_{\cz_r}} \\
\dTTo<{\alpha_!} && (6) & \rdTTo(4,2)^{\alpha_!} &&& \dTTo>{\alpha_!} \\
\dec{pt}{P_\cz} &&& \rTTo^{\iota^*} &&& \dec{pt}{\prod_rP_{\cz_r}}
\end{diagram}%
}
\caption{Coherence isomorphism $\coher$}
\label{dia-fig-coher-1-6}
\end{figure}

Coherence isomorphism in \figref{dia-fig-coher-1-6} is composed of
several isomorphisms. Isomorphism (1) exists by \propref{pro-pi-p-i}.
Isomorphism (2) follows from the fact that $q_*$ is an equivalence.

Isomorphism $T^{-2A}\to\rho_!\rho^*$ marked by (3) is obtained from
a sequence of affine linear $\prod P_{\cw_r}$-bundles
\begin{multline*}
U_\cw/(U_\cw\cap P_\cv) \rEpi \dots \rEpi
U_\cw/[U_\cw^{(k+1)}(U_\cw\cap P_\cv)] \\
\rEpi U_\cw/[U_\cw^{(k)}(U_\cw\cap P_\cv)] \rEpi \dots \rEpi pt,
\end{multline*}
where $U=U^{(1)}\supset U^{(2)}\supset\dots\supset U^{(m)}=1$ is
the lower central series of $U$. The group $\prod P_{\cw_r}$
acts by conjugation. Pull-back of this sequence along
$\prod\co''_r\times F_{\cw_r}|\prod P_{\cz_r}\times P_{\cw_r}\to
pt|\prod P_{\cw_r}$
gives a sequence of affine linear
$\prod P_{\cz_r}\times P_{\cw_r}$-bundles
\begin{multline*}
\prod\co''_r\times U_\cw/(U_\cw\cap P_\cv) \rEpi \dots \rEpi
\prod\co''_r\times U_\cw/[U_\cw^{(k+1)}(U_\cw\cap P_\cv)] \\
\rEpi \prod\co''_r\times U_\cw/[U_\cw^{(k)}(U_\cw\cap P_\cv)]
\rEpi \dots \rEpi \prod\co''_r.
\end{multline*}
Multiplication map $\prod\co''_r\times U_\cw\to\underline{\co}$ is an
isomorphism of $\prod P_{\cz_r}\times P_{\cw_r}$-spaces, hence, we may
replace the first space with the other one. The composition of the above
maps equals $\rho:\co/(U_\cw\cap P_\cv)\to\prod_r\co''_r$. The type (3)
isomorphism for $\rho$ is a composition of isomorphisms, constructed for
each affine bundle of the sequence in \corref{cor-aff-bun-equi}.

Isomorphism (4) is obtained from a commutative square of equivariant
maps. Isomorphism (6) is a base change isomorphism. Isomorphism (5) is
obtained from a base change isomorphism by inverting $q^*$ and $\pi^*$
and replacing them with their quasi-inverses $q_*$ and $\pi_*$.

We use the facts
\begin{gather*}
U_\cw = 1 + \oplus_{m,n;r>s} \Hom(V_r^m,V_s^n) ,\\
U_\cw\cap P_\cv = 1 + \oplus_{m\ge n;r>s} \Hom(V_r^m,V_s^n) ,\\
U_\cw/(U_\cw\cap P_\cv) \simeq \oplus_{m<n;r>s} \Hom(V_r^m,V_s^n),
\end{gather*}
The dimension of fiber bundle $\rho$
\begin{equation*}
A \overset{\text{def}}= \dim_\CC U_\cw/(U_\cw\cap P_\cv)
= \sum_{m<n;r>s} \dim V_r^m \cdot \dim V_s^n
\end{equation*}
multiplied by $-2$ equals to the total shift, obtained from all
braidings in the source functor of coherence isomorphism $\Coher$.

\subsection{Distinguished triangles}
Consider left $P_\cw$\n-invariant and right $P_\cv$\n-invariant subsets
$F\subset Y\subset G_V$, such that $F$ is closed in $Y$.
Denote $S=Y-F$ and consider inclusions
\begin{gather*}
F/{U_\cv} \rTTo^{i_U} Y/{U_\cv} \lTTo^{j_U} S/{U_\cv} ,\\
F/{P_\cv} \rTTo^{i_P} Y/{P_\cv} \lTTo^{j_P} S/{P_\cv}.
\end{gather*}
Denote $\phi^X$, $\pi^X$, $\alpha^X$ maps~\eqref{eq-p1p2p32} for
$\co=X\in\{Y,F,S\}$. Let $K\in D_{L_\cv}(pt)$,
$K_U=\phi^{Y*}K\in D_{P_\cw\times L_\cv}(Y/{U_\cv})$,
$K_P=\pi_{*}^YK_U\in D_{P_\cw}(Y/{P_\cv})$
(see the middle row in the following diagram).
\begin{diagram}
D_{L_\cv}(pt) & \rTTo^{\phi^{S*}} &
D_{P_\cw\times L_\cv}(S /{U_\cv}) &
\rTTo^{\pi_{*}^S}_\sim & D_{P_\cw}(S /{P_\cv}) &
\rTTo^{\alpha_{!}^S} & D_{P_\cw}(pt) \\
\dEq && \uTTo<{j_U^*} && \uTTo<{j_P^*} \dTTo>{j_{P!}} && \dEq \\
D_{L_\cv}(pt) & \rTTo^{\phi^{Y*}} &
D_{P_\cw\times L_\cv}(Y /{U_\cv}) &
\rTTo^{\pi_{*}^Y}_\sim & D_{P_\cw}(Y /{P_\cv}) &
\rTTo^{\alpha_{!}^Y} & D_{P_\cw}(pt) \\
\dEq && \dTTo<{i_U^*} && \dTTo<{i_P^*} \uTTo>{i_{P!}} && \dEq \\
D_{L_\cv}(pt) & \rTTo^{\phi^{F*}} &
D_{P_\cw\times L_\cv}(F /{U_\cv}) & \rTTo^{\pi_{*}^F}_\sim &
D_{P_\cw}(F /{P_\cv}) & \rTTo^{\alpha_{!}^F} & D_{P_\cw}(pt)
\end{diagram}
Applying $\alpha_{!}^Y$ to a standard triangle for $K_P$ we get a
distinguished triangle
\[ \alpha_{!}^Yj_{P!}j_P^*K_P \to \alpha_{!}^YK_P
\to \alpha_{!}^Yi_{P!}i_P^*K_P \to \]
The above diagram shows that it is isomorphic to certain following
sequences:
\[ \alpha_{!}^Sp_{2*}^Sj_U^*K_U \to \alpha_{!}^Y\pi_{*}^YK_U \to
\alpha_{!}^F\pi_{*}^Fi_U^*K_U \to \]
\[ \alpha_{!}^S\pi_{*}^S\phi^{S*}K \to \alpha_{!}^Y\pi_{*}^Y\phi^{Y*}K \to
\alpha_{!}^F\pi_{*}^F\phi^{F*}K \to \]
This is the triangle
\[ \Psi_\cw^{S;\cv}K \to \Psi_\cw^{Y;\cv}K \to
\Psi_\cw^{F;\cv}K \to \]
from which the structure triangle
\begin{equation}
\zhe_\cw^{S;\cv}K \rTTo^{a^{SY}} \zhe_\cw^{Y;\cv}K \rTTo^{b^{YF}}
\zhe_\cw^{F;\cv}K \rTTo^{c^{FS}}
\label{eq-trian-abc}
\end{equation}
is obtained by application of the exact functor $\pitchfork_\cw$.

This triangle replaces non-existing isomorphism
\(\zhe_\cw^{Y;\cv}K \simeq \zhe_\cw^{S;\cv}K \oplus
\zhe_\cw^{F;\cv}K \)
(this would be too much to ask for).

\begin{proposition}
For any pair of closed embeddings $F\subset Z\subset Y$ set
$S=Y-F$, $Q=Z-F$, $R=Y-Z$. The following diagram made with
distinguished triangles \eqref{eq-trian-abc} is an octahedron.
\begin{equation*}
\begin{diagram}[inline,width=2.8em]
\zhe^{R} && \rTTo^{a^{RY}} && \zhe^{Y} \\
& \luTTo^1_{c^{ZR}} & d & \ldTTo_{b^{WZ}} & \\
\uTTo<{c^{QR}}>1 & = & \zhe^{Z} & = & \dTTo>{b^{YF}} \\
& \ruTTo^{a^{QZ}} & d & \rdTTo^{b^{ZF}} & \\
\zhe^{Q} && \lTTo_1^{c^{FQ}} && \zhe^{F}
\end{diagram}
\quad
\begin{diagram}[inline,width=2.8em]
\zhe^{R} && \rTTo^{a^{RY}} && \zhe^{Y} \\
& \rdTTo_{a^{RS}} & = & \ruTTo_{a^{SY}} & \\
\uTTo<{c^{QR}}>1 & d & \zhe^{S} & d & \dTTo>{b^{YF}} \\
& \ldTTo^{b^{SQ}} & = & \luTTo_1^{c^{FS}} & \\
\zhe^{Q} && \lTTo_1^{c^{FQ}} && \zhe^{F}
\end{diagram}
\label{dia-zhe-octa}
\end{equation*}
That is, the four triangles marked ``='' commute,
as well as two squares
\[
\begin{diagram}[inline]
\zhe^Z & \rTTo^{b^{ZF}} & \zhe^F \\
\dTTo<{c^{ZR}} && \dTTo>{c^{FS}} \\
T\zhe^R & \rTTo^{Ta^{RS}} & T\zhe^S
\end{diagram}
\qquad
\begin{diagram}[inline]
\zhe^S & \rTTo^{a^{SY}} & \zhe^Y \\
\dTTo<{b^{SQ}} && \dTTo>{b^{YZ}} \\
\zhe^Q & \rTTo^{a^{QZ}} & \zhe^Z
\end{diagram}
\quad.
\]
\end{proposition}

The proof follows from \corref{cor-octa-equiv}.

\appendix
\section{Technical results}
\subsection{An equivalence of equivariant derived categories}
\label{sec-equi-equi}
Let ${P}$ be a linear complex algebraic group, let $U$, ${L}$ be its
subgroups such that ${U}$ is normal in ${P}$ and the manifold $U$ is
an affine space. Assume that $P$ is a semidirect product of $U$ and $L$,
that is, the multiplication map $\cdot:U\times L\to P$ is an isomorphism
of manifolds. Denote by $\iota:L\rMono P$ and $\kappa:P\rEpi L$
the natural inclusion and projection.
Let ${Y}$ be a ${P}$\n-variety, let ${E}$ be an
${L}$\n-variety, let $i:{E}\to{Y}$ be a
${\iota}$\n-map  and let $p:{Y}\to{E}$ be a
${\kappa}$\n-map identifying ${E}$ with the quotient
${U}\backslash{Y}$ such that $p\circ i=\id_E$. Consider
\[ W = {P}\times_{{L}}{E} = ({P}\times{E})/{L}, \quad
(p,e).l \overset{\text{def}}= (pl,l^{-1}e) .\]
The variety $W$ is a ${P}$\n-variety via left translations on ${P}$.
As a ${P}$\n-variety $W$ is isomorphic to ${U}\times{E}$,
where ${P}$ acts by
\[ (u,l).(v,e) = (u\cdot(lvl^{-1}),l.e) .\]
There is a surjective map of ${P}$\n-varieties
\[ \pi:W\simeq{U}\times{E}\to{Y},
\quad (u,e) \mapsto u.i(e) .\]
Assume furthermore that $\pi$ is a locally trivial fibration with
affine fibers. Then $\pi:W\to{Y}$ is $\infty$\n-acyclic in the sense of
\cite{BernsL:Equivariant}.

\begin{lemma}\label{lem-pi-equiv}
In the above assumptions the functor
$\pi^*:D_{{P}}^b({Y})\to D_{{P}}^b(W)$ is fully faithful.
\end{lemma}

\begin{proof}
Choose an interval $I=[a,b]\subset\ZZ$ and a number $m\geq|I|$.
Let $M=M_m$ be an $m$\n-acyclic free ${P}$\n-manifold. Then
$R=M\times W\rTTo^{\pr}W$ is a $m$\n-acyclic ${P}$\n-resolution
of $W$ and $r:R \rTTo^{\pr_2} W \rTTo^\pi Y$ is such of
${Y}$. We know that
\[ \pi^*: D^I({Y}) \to D^I(W) \]
is full and faithful, see Proposition 1.9.2 of
\cite{BernsL:Equivariant}. Therefore, functor from
Section~2.1.6 of \cite{BernsL:Equivariant}
\begin{equation}
\pi^*: D_{{P}}^I({Y},R) \to D_{{P}}^I(W,R),
\label{eq-pi-res}
\end{equation}
\begin{multline*}
\pi^*(F_{{Y}},\overline{F},
\beta:r^*F_{{Y}}\rTTo^\sim q^*\overline{F}) \\
= (\pi^*F_{{Y}},\overline{F},\gamma:\pr^*\pi^*F_{{Y}}\simeq
r^*F_{{Y}} \rTTo_\sim^\beta q^*\overline{F}),
\end{multline*}
is also full and faithful, where $q:R\to\overline{R}$,
$r:R\rTTo^\pr W\rTTo^\pi Y$.
Indeed, a morphism of $D_{{P}}^I({Y},R)$
is a pair $(\alpha_{{Y}}:F_{{Y}}\to
H_{{Y}},\overline{\alpha}:\overline{F}\to\overline{H})$,
which makes commutative the right square in the following diagram
\begin{diagram}[LaTeXeqno]
\pr^*\pi^*F_{{Y}} & \rTTo^\sim & r^*F_{{Y}} & \rTTo^\beta &
q^*\overline{F} \\
\dTTo<{\pr^*\pi^*\alpha_{{Y}}} & = &
\dTTo<{r^*\alpha_{{Y}}} && \dTTo>{q^*\overline{\alpha}} & . \\
\pr^*\pi^*H_{{Y}} & \rTTo^\sim & r^*H_{{Y}} & \rTTo^\beta &
q^*\overline{H}
\label{dia-mor-D(W,R)}
\end{diagram}
The morphism $(\alpha_{{Y}},\overline{\alpha})$ is sent by
$\pi^*$ to a morphism $(\pi^*\alpha_{{Y}},\overline{\alpha})$.
This shows faithfulness of functor~\eqref{eq-pi-res}. This functor is
also full because any morphism
$(\pi^*F_{{Y}},\overline{F},\phi) \to (\pi^*H_{{Y}},\overline{H},\eta)$
has to be of the form $(\pi^*\alpha_{{Y}},\overline{\alpha})$ for some
$\alpha_{{Y}}$. The commutativity of the exterior of
diagram~\eqref{dia-mor-D(W,R)} implies that the pair
$(\alpha_{{Y}},\overline{\alpha})$ is a morphism.

Assign $D_{{P}}^I({Y})=D_{{P}}^I({Y},R)$,
$D_{{P}}^I(W)=D_{{P}}^I(W,R)$ and vary $R=M_m\times W$, increasing $m$.
Since $D_{{P}}^b({Y})$ (resp. $D_{{P}}^b({W})$)
is a colimit of a system of fully faithful functors
$D_{{P}}^J({Y})\to D_{{P}}^I({Y})$, $J\subset I$,
(by Definition~2.2.4 of \cite{BernsL:Equivariant}),
we see that the colimits are also fully and faithfully embedded.
\end{proof}

\begin{proposition}\label{pro-pi-p-i}
The functors
\begin{align*}
\pi^*:& \Dbc_{{P}}({Y}) \to \Dbc_{{P}}(W), \\
p^*:& \Dbc_{{L}}({E}) \to \Dbc_{{P}}({Y}), \\
i^*:& \Dbc_{{P}}({Y}) \to \Dbc_{{L}}({E})
\end{align*}
are equivalences, and the last two are quasi-inverse to each other.
\end{proposition}

\begin{proof}
Notice that $W\simeq {U}\times{E}$ is a free ${U}$\n-variety.
Denote by $\phi:{P} \rEpi L={P}/{U}$ the canonical projection.
The $\phi$\n-map $s:W\rTTo^\pi Y\rTTo^p E$ is
the projection $W\to{U}\backslash W={E}$.

Therefore, the functor $s^*:\Dbc_{{L}}({E})\to \Dbc_{{P}}(W)$ is an
equivalence by Theorem~2.6.2 of \cite{BernsL:Equivariant}.
It decomposes as
\[ s^* = \bigl( \Dbc_{{L}}({E}) \rTTo^{p^*}
\Dbc_{{P}}({Y}) \rTTo^{\pi^*} \Dbc_{{P}}(W) \bigr) .\]
We conclude that $\pi^*$ is essentially surjective on objects, hence,
an equivalence by \lemref{lem-pi-equiv}. Therefore, $p^*$ is also an
equivalence. Since
\[ \bigl( \Dbc_{{L}}({E}) \rTTo^{p^*} \Dbc_{{P}}({Y})
\rTTo^{i^*} \Dbc_{{L}}({E}) \bigr) \simeq \Id \]
$i^*$ is a quasi-inverse of $p^*$ (and also an equivalence).
\end{proof}

\section{Properties of standard isomorphisms}
\subsection{Inverse image isomorphisms}
\begin{proposition}\label{pro-iota-cocycle}
Isomorphisms $\iota$ from \eqref{eq-dia-iota} satisfy
the cocycle condition:
\begin{equation}
\begin{diagram}[inline,width=2.6em,tight]
\Dbc_K(Z) && \rTTo^{g^*} && \Dbc_H(Y) \\
& \rdTwoar_{\iota_{h,g}} && \ruTTo(4,4)_{(hg)^*}
\ldTwoar(2,4)_{\iota_{hg,f}} & \\
\uTTo<{h^*} &&&& \dTTo>{f^*} \\
&&&& \\
\Dbc_L(W) && \rTTo_{(hgf)^*} && \Dbc_G(X)
\end{diagram}
\quad=\quad
\begin{diagram}[inline,width=2.6em,tight]
\Dbc_K(Z) && \rTTo^{g^*} && \Dbc_H(Y) \\
& \rdTwoar(2,4)_{\iota_{h,gf}} \rdTTo(4,4)_{(gf)^*} &&
\ldTwoar_{\iota_{g,f}} & \\
\uTTo<{h^*} &&&& \dTTo>{f^*} \\
&&&& \\
\Dbc_L(W) && \rTTo_{(hgf)^*} && \Dbc_G(X)
\end{diagram}
\label{eq-cocyc-iota}
\end{equation}
\end{proposition}

\begin{proof}
Compose this equation with the functor
\[\Dbc_G(X)\to\Dbc(\Resto(f,g,h)\to\Top)\]
and substitute the definition of
$\iota$. Its interpretation might be the following: $\iota$ and $i$ are
gauge equivalent via gauge transformation $\eta$. The required equation
reduces to equality of the following two isomorphisms:
\begin{diagram}[width=5.5em,tight,nobalance]
&& \Dbc(\Resto(h)\to\Top) & \rTTo^{g^*} & \Dbc(\Resto(g,h)\to\Top) \\
& \ruTTo^{h^*} & \dTwoar>{i_{h,g}} & \ruTTo_\sim & \dTTo>{f^*} \\
\Dbc_L(W) & \rTTo^{(hg)^*} & \Dbc(\Resto(hg)\to\Top) & = &
\Dbc(\Resto(f,g,h)\to\Top) \\
\dTTo<{(hgf)^*} & \ldTwoar_{i_{hg,f}} & \dTTo>{f^*} & \ruTTo_\sim & \\
\Dbc(\Resto(hgf)\to\Top) & \rTTo_\sim & \Dbc(\Resto(f,hg)\to\Top) &&
\end{diagram}
\begin{diagram}[width=5.5em,tight,nobalance]
&& \Dbc(\Resto(h)\to\Top) & \rTTo^{g^*} & \Dbc(\Resto(g,h)\to\Top) \\
& \ruTTo^{h^*} \ldTwoar(2,4)^{i_{h,gf}} & \dTTo>{(gf)^*} &
\ldTwoar_{i_{g,f}} & \dTTo>{f^*} \\
\Dbc_L(W) && \Dbc(\Resto(gf,h)\to\Top) & \rTTo_\sim &
\Dbc(\Resto(f,g,h)\to\Top) \\
\dTTo<{(hgf)^*} & \ruTTo_\sim & = & \ruTTo_\sim & \\
\Dbc(\Resto(hgf)\to\Top) & \rTTo_\sim & \Dbc(\Resto(f,hg)\to\Top) &&
\end{diagram}
This equation follows from \eqref{eq-cocyc-iota} for $\Dbc(T)$,
$T\in\Top$, which in turn follows from the same equation for $D^b(T)$
and, in turn,  for the category of sheaves.
\end{proof}

\subsection{Properties of base change isomorphisms}
\begin{lemma}\label{lem-distr-base-chge}
Whenever squares 1 and 2 below are pull-back squares
of algebraic varieties
\begin{diagram}[LaTeXeqno]
X\SEpbk & \rTTo^\kappa & Y\SEpbk & \rTTo^\lambda & Z \\
\dTTo<f & 1 & \dTTo<g & 2 & \dTTo>h \\
U & \rTTo^\nu & V & \rTTo^\mu & W
\label{dia-2-pbk-squares}
\end{diagram}
the following equation holds
\begin{diagram}[LaTeXeqno]
\nu^*\mu^*(Rh_!) & \rTTo^{\nu^*\beta_2} & \nu^*(Rg_!)\lambda^* &
\rTTo^{\beta_1} & (Rf_!)\kappa^*\lambda^* \\
\dTTo<{i_{\mu,\nu}} && = && \dTTo>{Rf_!i_{\lambda,\kappa}} && . \\
(\mu\nu)^*Rh_! && \rTTo^\beta && (Rf_!)(\lambda\kappa)^*
\label{dia-eq-4-bas-chg}
\end{diagram}
The statement can be reformulated as an equation
\begin{equation*}
\begin{diagram}[inline,width=2.6em,tight]
\Dbc(X) && \lTTo^{(\lambda\kappa)^*} && \Dbc(Z) \\
&&&& \\
\dTTo<{Rf_!} &&&& \dTTo>{Rh_!} \\
\Dbc(U) && \lTTo^{(\mu\nu)^*} && \Dbc(W) \\
& \luTTo_{\nu^*} && \ldTTo_{\mu^*} & \\
&& \Dbc(V) &&
\end{diagram}
\quad=\quad
\begin{diagram}[inline,width=2.6em,tight]
\Dbc(X) && \lTTo^{(\lambda\kappa)^*} && \Dbc(Z) \\
& \luTTo_{\kappa^*} && \ldTTo_{\lambda^*} & \\
\dTTo<{Rf_!} && \Dbc(Y) && \dTTo>{Rh_!} \\
\Dbc(U) && \dTTo<{Rg_!} && \Dbc(W) \\
& \luTTo_{\nu^*} && \ldTTo_{\mu^*} & \\
&& \Dbc(V) &&
\end{diagram}
\end{equation*}
\end{lemma}

\begin{proof}
First we prove analogous to \eqref{dia-eq-4-bas-chg} equation for
sheaves. We use the definition of the base change isomorphism
of sheaves, corresponding to the pull-back square
\begin{diagram}
\0\SEpbk & \rTTo^e & \0 \\
\dTTo<h && \dTTo>l \\
\0 & \rTTo_k & \0
\end{diagram}
It is defined as the composition
\[ \beta: k^*l_! \rTTo^\eta k^*l_!e_*e^* \rMono
k^*k_* h_!e^* \rTTo^\epsilon h_!e^* \]
or as
\[ \beta =\quad
\begin{diagram}[inline]
\0 & \lTTo^{e^*} & \0 \\
\dTTo<{h_!} & \rdTTo_{e_*} & \dId<{\overset\eta\Longleftarrow} \\
\0 & \lTwoar & \0 \\
\dId>{\overset\epsilon\Longleftarrow} & \rdTTo^{k_*} & \dTTo>{l_!} \\
\0 & \lTTo_{k^*} & \0
\end{diagram}
\quad,
\]
where $\epsilon$ and $\eta$ are adjunction maps.
Equation~\eqref{dia-eq-4-bas-chg} for sheaves which we want to prove
becomes an equation between the following two isomorphisms
\begin{diagram}[LaTeXeqno]
\Sh X && \lTTo^{(\lambda\kappa)^*} && \Sh Z \\
\dTTo<{f_!} & \rdTTo(4,2)_{(\lambda\kappa)_*} &&&
\dId<{\overset\eta\Longleftarrow} \\
\Sh U && \lTwoar && \Sh Z \\
\dId>{\overset\epsilon\Longleftarrow} & \rdTTo(4,2)^{(\mu\nu)_*} &&&
\dTTo>{h_!} \\
\Sh U && \lTTo^{(\mu\nu)^*} && \Sh W \\
& \luTTo_{\nu^*} & \uTwoar & \ldTTo_{\mu^*} & \\
&& \Sh V &&
\label{dia-Sh-ustar-lshriek}
\end{diagram}
\[
\begin{diagram}[inline,width=2.75em,tight]
\Sh X && \lTTo^{(\lambda\kappa)^*} && \Sh Z \\
& \rdTTo(2,4)_{\kappa_*} \luTTo^{\kappa^*} & \uTwoar &
\ldTTo^{\lambda^*} & \\
\dTTo<{f_!} && \Sh Y && \dId<{\overset\eta\Longleftarrow} \\
&& \dId<{\overset\Longleftarrow\eta} && \\
\Sh U & \lTwoar & \Sh Y & \rTTo^{\lambda_*} & \Sh Z \\
\dId>{\overset\epsilon\Longleftarrow} & \rdTTo_{\nu_*} &
\dTTo>{g_!} & \ldTwoar & \dTTo>{h_!} \\
\Sh U && \Sh V & \rTTo_{\mu_*} & \Sh W \\
& \luTTo_{\nu^*} & \dId>{\overset\epsilon\Longleftarrow} &
\ldTTo_{\mu^*} & \\
&& \Sh V &&
\end{diagram}
\quad=\quad
\begin{diagram}[inline,width=2.75em,tight]
\Sh X && \lTTo^{(\lambda\kappa)^*} && \Sh Z \\
& \rdTTo_{\kappa_*} \rdTTo(4,2)^{(\lambda\kappa)_*} &&&
\dId<{\overset\eta\Longleftarrow}  \\
\dTTo<{f_!} && \Sh Y & \rTTo_{\lambda_*} & \Sh Z \\
& \ldTwoar & \dTTo>{g_!} & \ldTwoar & \\
\Sh U & \rTTo^{\nu_*} & \Sh V && \dTTo>{h_!} \\
\dId>{\overset\epsilon\Longleftarrow} &
\rdTTo(4,2)[nohug]_{(\mu\nu)_*}^= && \rdTTo^{\mu_*} & \\
\Sh U && \lTTo_{(\mu\nu)^*} && \Sh W \\
& \luTTo_{\nu^*} & \uTwoar & \ldTTo_{\mu^*} & \\
&& \Sh V &&
\end{diagram}
\]
The parallelogram in diagram \eqref{dia-Sh-ustar-lshriek} equals to the
parallelogram in the right diagram above since the embedding
$h_!(\lambda\kappa)_* \rMono (\mu\nu)_*f_!$ can be presented as a
composition of two embeddings
\[ h_!\lambda_*\kappa_* \rMono \mu_*g_!\kappa_* \rMono \mu_*\nu_*f_!. \]
Indeed, all these functors are subfunctors of
$h_*\lambda_*\kappa_* = \mu_*\nu_*f_*$.

        From equation \eqref{dia-eq-4-bas-chg} for sheaves
we deduce it for derived functors.
\end{proof}

\begin{proposition}\label{prop-distr-base-chge}
Let $f:X\to U$, $g:Y\to V$ and $h:Z\to W$ be $G$\n-equivariant,
$H$\n-equivariant and $K$\n-equivariant maps of algebraic varieties.
Let $X\rTTo^\kappa Y \rTTo^\lambda Z$ and $U \rTTo^\nu V \rTTo^\mu W$
be $G \rTTo^\phi H \rTTo^\psi K$\n-equivariant.
Assume that squares 1 and 2 of commutative diagram
\eqref{dia-2-pbk-squares} are pull-back squares.
Then the following equation holds
\begin{diagram}
\nu^*\mu^*h_! & \rTTo^{\nu^*\beta_2} & \nu^*g_!\lambda^* &
\rTTo^{\beta_1} & f_!\kappa^*\lambda^* \\
\dTTo<{i_{\mu,\nu}} && = && \dTTo>{f_!i_{\lambda,\kappa}} \\
(\mu\nu)^*h_! && \rTTo^\beta && f_!(\lambda\kappa)^*
\end{diagram}
The statement can be reformulated as an equation
\begin{equation*}
\begin{diagram}[inline,width=2.6em,tight]
\Dbc_G(X) && \lTTo^{(\lambda\kappa)^*} && \Dbc_K(Z) \\
&&&& \\
\dTTo<{f_!} &&&& \dTTo>{h_!} \\
\Dbc_G(U) && \lTTo^{(\mu\nu)^*} && \Dbc_K(W) \\
& \luTTo_{\nu^*} && \ldTTo_{\mu^*} & \\
&& \Dbc_H(V) &&
\end{diagram}
\quad=\quad
\begin{diagram}[inline,width=2.6em,tight]
\Dbc_G(X) && \lTTo^{(\lambda\kappa)^*} && \Dbc_K(Z) \\
& \luTTo_{\kappa^*} && \ldTTo_{\lambda^*} & \\
\dTTo<{f_!} && \Dbc_H(Y) && \dTTo>{h_!} \\
\Dbc_G(U) && \dTTo<{g_!} && \Dbc_K(W) \\
& \luTTo_{\nu^*} && \ldTTo_{\mu^*} & \\
&& \Dbc_H(V) &&
\end{diagram}
\end{equation*}
\end{proposition}

\begin{proof}
The statement reduces to equality between the following two
isomorphisms:
\begin{diagram}[width=4em,nobalance]
\Dbc(\Resto(\kappa,\lambda;\phi,\psi)\to\Top) & \lTTo^\sim &
\Dbc(\Resto(\lambda\kappa;\psi\phi)\to\Top) & \lTTo^{(\lambda\kappa)^*}
& \Dbc_K(Z) \\
\dTTo<{f_!} & = & \dTTo<{f_!} & \luTwoar^\beta & \dTTo>{h_!} \\
\Dbc(\Resto(\nu,\mu;\phi,\psi)\to\Top) & \lTTo^\sim &
\Dbc(\Resto(\mu\nu;\psi\phi)\to\Top) & \lTTo^{(\mu\nu)^*}
& \Dbc_K(W) \\
& \luTTo_{\nu^*} & \uTwoar>i & \ldTTo_{\mu^*} & \\
&& \Dbc(\Resto(\mu,\psi)\to\Top) &&
\end{diagram}
\begin{diagram}[width=4em,nobalance]
\Dbc(\Resto(\kappa,\lambda;\phi,\psi)\to\Top) & \lTTo^\sim &
\Dbc(\Resto(\lambda\kappa;\psi\phi)\to\Top) & \lTTo^{(\lambda\kappa)^*}
& \Dbc_K(Z) \\
\dTTo<{f_!} & \luTwoar(2,4)>\beta \luTTo_{\kappa^*} & \uTwoar>i &
\ldTTo_{\lambda^*} & \dTTo>{h_!} \\
\Dbc(\Resto(\nu,\mu;\phi,\psi)\to\Top) &&
\Dbc(\Resto(\lambda,\phi)\to\Top) & \lTwoar^\beta & \Dbc_K(W) \\
& \luTTo_{\nu^*} & \dTTo<{g_!} & \ldTTo_{\mu^*} & \\
&& \Dbc(\Resto(\mu,\psi)\to\Top) &&
\end{diagram}
Substituting definitions of isomorphisms $\beta$ from
\secref{sec-equi-bas-chg} and similar, we see that the
equation between the above isomorphisms follows from
\lemref{lem-distr-base-chge}, that is,
from the non-equivariant version.
\end{proof}

\section{Distinguished triangles and octahedrons}
\subsection{Standard distinguished triangles.}
\label{equi-dis-tria}
Let us consider a closed $G$\n-invariant subspace $i:X\rMono Y$ with
the complement $S=Y-X$, $j:S\rMono Y$. Let $P\to Y$ be a resolution of
$Y$. As explained in \cite{BernsL:Equivariant} resolutions of $X$
of the form $X\times_YP\to X$ taken from the pull-back square
\begin{diagram}
X\times_YP\SEpbk & \rTTo^{\underline{i}} & X \\
\dTTo && \dTTo \\
X & \rTTo^i & Y
\end{diagram}
suffice to determine an object of $\Dbc_G(X)$ up to an isomorphism.
Therefore, the functor $i_!i^*$ is canonically isomorphic to the functor
\begin{align*}
\widetilde{i_!i^*}: \Dbc_G(Y) & \longrightarrow \Dbc_G(Y) \\
{}[P\mapsto K(P)] & \longmapsto
[P\mapsto \overline{\underline{i}}_!\overline{\underline{i}}^*(K(P))]
\end{align*}
Similarly for $j_!j^*$.

There is a distinguished triangle in $\Dbc_G(Y)$
\cite{BernsL:Equivariant}
\[ \widetilde{j_!j^*}K \rTTo^{a} K \rTTo^{b}
\widetilde{i_!i^*}K \rTTo^{c} \]
given by a collection of standard distinguished triangles
\cite{BeilBD:Perverse}
\[ \overline{\underline{j}}_!\overline{\underline{j}}^*(KP)
\rTTo^{a_{KP}} KP \rTTo^{b_{KP}}
\overline{\underline{i}}_!\overline{\underline{i}}^*(KP)
\rTTo^{c_{KP}} \]
for $P\in\Resto(Y,G)$. The canonically isomorphic to it distinguished
triangle in $\Dbc_G(Y)$ is denoted
\[ j_!j^*K \rTTo^{a_K} K \rTTo^{b_K} i_!i^*K \rTTo^{c_K} \quad.\]

\subsection{Standard octahedron associated with a triple of
spaces.}
Suppose that $X\subset Z\subset Y$ are closed embeddings of
stratified subspaces of some stratified space and $ S= Y- X$,
$ Q= Z- X$, $ R= Y- Z$. With a closed embedding $i$ and the
complementary embedding $j$ is associated a standard triangle
\cite{BeilBD:Perverse} in $D(X)$
\[ j_!j^*K \rTTo^{a_K} K \rTTo^{b_K} i_!i^*K \rTTo^{c_K}\quad.
\]
Denote by $i_{YZ}$, $i_{YX}$, $i_{ZX}$, $i_{SQ}$ the closed embeddings
and $j_{RY}$, $j_{QZ}$, $j_{SY}$, $j_{RY}$ the open ones. Add
corresponding indices to the triangle maps $a$, $b$, $c$.

\begin{proposition}\label{pro-stan-octa}
The obvious identification of vertices makes the following diagram
in $D(Y)$ into an octahedron
\begin{diagram}[width=5.4em]
j_{RY!}j^*_{RY}K && \rTTo^{a^{RY}} && K \\
& \luTTo_{c^{ZR}} && \ldTTo_{b^{YZ}} & \\
\uTTo<1 & \hspace*{-5mm} = \hspace*{5mm} & i_{YZ!}i_{YZ}^*K &
\hspace*{5mm} = \hspace*{-5mm} & \dTTo \\
& \ruTTo^{i_{YZ!}a^{QZ}} && \rdTTo^{i_{YZ!}b^{ZX}} & \\
i_{YZ!}j_{QZ!}j_{QZ}^*i_{YZ}^*K && \lTTo^{i_{YZ!}c^{XQ}} &&
i_{YZ!}i_{ZX!}i_{ZX}^*i_{YZ}^*K
\end{diagram}
\begin{diagram}[width=5.4em]
j_{SY!}j_{RS!}j_{RS}^*j_{SY}^*K && \rTTo && K \\
& \rdTTo^{j_{SY!}a^{RS}} & = & \ruTTo^{a^{SY}} & \\
\uTTo<{j_{SY!}c^{QR}} && j_{SY!}j_{SY}^*K && \dTTo>{b^{YX}} \\
& \ldTTo_{j_{SY!}b^{SQ}} & = & \luTTo_{c^{XS}} & \\
j_{SY!}i_{SQ!}i_{SQ}^*j_{SY}^*K && \lTTo_1 && i_{YX!}i_{YX}^*K
\end{diagram}
That is, the following ``triangles'' and ``squares'' commute.
\begin{diagram}
Tj_{SY!}j_{RS!}j_{RS}^*j_{SY}^*K & \rTTo_\sim^{\iota\cdot\iota} &
Tj_{RY!}j_{RY}^*K && \\
&&& \luTTo^{c^{ZR}} & \\
\uTTo<{j_{SY!}c^{QR}} && 1 && i_{YZ!}i_{YZ}^*K \\
&&& \ruTTo_{i_{YZ!}a^{QZ}} & \\
j_{SY!}j_{SQ!}j_{SQ}^*j_{SY}^*K & \rTTo_\sim &
i_{YZ!}j_{QZ!}j_{QZ}^*i_{YZ}^*K &&
\end{diagram}
\begin{diagram}
i_{YZ!}i_{YZ}^*K & \lTTo^{b^{YZ}} & K \\
\dTTo<{i_{YZ!}b^{YZ}} & 2 & \dTTo>{b^{YX}} \\
i_{YZ!}i_{ZX!}i_{ZX}^*i_{YZ}^*K & \lTTo_\sim^{\iota\cdot\iota} &
i_{YX!}i_{YX}^*K
\end{diagram}
\begin{diagram}
j_{RY!}j_{RY}^*K & \rTTo^{a^{RY}} & K \\
\dTTo<{\iota\cdot\iota}>\wr & 3 & \uTTo>{a^{SY}} \\
j_{SY!}j_{RS!}j_{RS}^*j_{SY}^*K & \rTTo^{j_{SY!}a^{RS}} &
j_{SY!}j_{SY}^*K
\end{diagram}
\begin{diagram}
&& Tj_{SY!}j_{SY}^*K && \\
& \ldTTo^{j_{SY!}b^{SQ}} && \luTTo^{c^{XS}} & \\
Tj_{SY!}i_{SQ!}i_{SQ}^*j_{SY}^*K && 4 && i_{YX!}i_{YX}^*K \\
\dTTo<\wr &&&& \dTTo<\wr>{\iota\cdot\iota} \\
Ti_{YZ!}j_{QZ!}j_{QZ}^*i_{YZ}^*K && \lTTo^{i_{YZ!}c^{XQ}} &&
i_{YZ!}i_{ZX!}i_{ZX}^*i_{YZ}^*K
\end{diagram}
\begin{diagram}
i_{YZ!}i_{YZ}^*K & \rTTo^{i_{YZ!}b^{ZX}} &
i_{YZ!}i_{ZX!}i_{ZX}^*i_{YZ}^*K & \rTTo_\sim^{\iota\cdot\iota} &
i_{YX!}i_{YX}^*K \\
\dTTo<{c^{ZR}} && 5 && \dTTo>{c^{XS}} \\
Tj_{RY!}j_{RY}^*K & \rTTo_\sim^{\iota\cdot\iota} &
Tj_{SY!}j_{RS!}j_{RS}^*j_{SY}^*K & \rTTo^{j_{SY!}a^{RS}} &
Tj_{SY!}j_{SY}^*K
\end{diagram}
\begin{diagram}
j_{SY!}j_{SY}^*K && \rTTo^{a^{SY}} && K \\
\dTTo<{j_{SY!}b^{SQ}} && 6 && \dTTo>{b^{YZ}} \\
j_{SY!}i_{SQ!}i_{SQ}^*j_{SY}^*K & \rTTo_\sim &
i_{YZ!}j_{QZ!}j_{QZ}^*i_{YZ}^*K & \rTTo^{i_{YZ!}a^{QZ}} &
i_{YZ!}i_{YZ}^*K
\end{diagram}
\end{proposition}

\begin{proof}
This octahedron in a derived category of sheaves follows
\cite{BeilBD:Perverse} from the commutative diagram with exact lines in
the abelian category of sheaves.
\begin{diagram}[width=2em,height=1em,tight,LaTeXeqno]
&&&& 0 &&&& 0 &&&& \\
&&&&& \rdTTo(2,6) && \ruTTo(2,6) &&&&& \\
&&&& &&&& &&&& \\
&&&& &&&& &&&& \\
&&&& &&&& &&&& \\
&&&& &&&& &&&& \\
&&&&&& k_{QY!}k_{QY}^*\cf &&&&&& \\
&&&&& \ruTTo(2,6)^{j_{SY!}b^{SQ}} &&
\rdTTo(2,6)^{i_{YZ!}a^{QZ}} &&&&& \\
&&&& &&&& &&&& \\
&&&& &&&& &&&& \\
&& 0 &&&& &&&& 0 && \\
&&& \rdTTo &&&&&& \ruTTo &&& \\
&&&& j_{SY!}j_{SY}^*\cf &&&& i_{YZ!}i_{YZ}^*\cf &&&& \\
&&& \ruTTo(2,6)^{j_{SY!}a^{RS}} && \rdTTo_{a^{SY}} &&
\ruTTo_{b^{YZ}} && \rdTTo(2,6)^{i^{YZ!}b^{ZX}} &&& \\
&&&&&& \cf &&&&&& \\
&&&&& \ruTTo(4,4)_{a^{RY}} && \rdTTo(4,4)_{b^{YX}} &&&&& \\
&&&& &&&& &&&& \\
&&&& &&&& &&&& \\
&& j_{RY!}j_{RY}^*\cf &&&& &&&& i_{YX!}i_{YX}^*\cf && \\
& \ruTTo \ruTTo(2,6) &&&&& &&&&& \rdTTo(2,6) \rdTTo & \\
0 &&&& &&&& &&&& 0 \\
&&&& &&&& &&&& \\
&&&& &&&& &&&& \\
&&&& &&&& &&&& \\
0 &&&& &&&& &&&& 0
\label{dia-wigwam}
\end{diagram}
Here $k_{QY}: Q\rMono Y$ denotes the embedding. The ambiguous
north-north-east line is isomorphic to the exact sequence
\[ 0 \to j_{SY!}j_{RS!}j_{RS}^*j_{SY}^*\cf
\rTTo^{j_{SY!}a^{RS}_{j^*\cf}} j_{SY!}j_{SY}^*\cf
\rTTo^{j_{SY!}b^{SQ}_{j^*\cf}} j_{SY!}i_{SQ!}i_{SQ}^*j_{SY}^*\cf
\to 0 .\]
The ambiguous south-south-east line is isomorphic to the exact sequence
\[ 0 \to i_{YZ!}j_{QZ!}j_{QZ}^*i_{YZ}^*\cf
\rTTo^{i_{YZ!}a^{QZ}_{i^*\cf}} i_{YZ!}i_{YZ}^*\cf
\rTTo^{i_{YZ!}b^{ZX}_{i^*\cf}} i_{YZ!}i_{ZX!}i_{ZX}^*i_{YZ}^*\cf
\to 0 .\]
This shows exactness of the lines of the diagram. It remains to
show commutativity of the two triangles and the quadrangle.

Together with a hidden isomorphism the left triangle is square 3:
\begin{diagram}
j_{SY!}j_{RS!}j_{RS}^*j_{SY}^*\cf & \rTTo^{j_{SY!}a_{j^*\cf}^R} &
j_{SY!}j_{SY}^*\cf \\
\dTTo<\wr && \dTTo>{a_\cf^{SY}} \\
j_{RY!}j_{RY}^*\cf & \rTTo^{a_\cf^{RY}} & \cf
\label{dia-left-triangle}
\end{diagram}
Recall that for any sheaf $\cf$ on $Y=S\cup X$ we have
\( j_{SY!}j_{SY}^*\cf = \cf_S = (\cf\big|_S)^Y \) and
\( i_{YX!}i_{YX}^*\cf = \cf_X = (\cf\big|_X)^Y \), the sheaves extended
by $0$ from $S$ (resp. from $X$). The embedding
$a_\cf^{SY}:j_{SY!}j_{SY}^*\cf\to\cf$ is the embedding of
\'espaces \'etal\'es
$\cf_S\rMono\cf$. Therefore, \eqref{dia-left-triangle} states simply
that the triangle of embeddings
\begin{diagram}
&& \cf_S && \\
& \ruMono && \rdMono & \\
\cf_R && \rMono && \cf
\end{diagram}
is commutative, which is obvious.

Similarly, the right triangle of diagram~\eqref{dia-wigwam} is square 4.
\begin{diagram}
\cf & \rTTo^{b_\cf^{YZ}} & i_{YZ!}i_{YZ}^*\cf \\
\dTTo<{b_\cf^{YX}} && \dTTo>{i_{YZ!}b_{i^*\cf}^{ZX}} \\
i_{YX!}i_{YX}^*\cf & \rTTo^\sim & i_{YZ!}i_{ZX!}i_{ZX}^*i_{YZ}^*\cf
\end{diagram}
Since the surjection $b_\cf^{YX}:\cf\to i_{YX!}i_{YX}^*\cf$ is
interpreted as the projection $\cf\rEpi\cf_X$ of \'espaces \'etal\'es,
the commutativity of this square reduces to commutativity of
\begin{diagram}
&& \cf_Z && \\
& \ruEpi && \rdEpi & \\
\cf && \rEpi && \cf_X
\end{diagram}
which is obvious.

The quadrangle in diagram~\eqref{dia-wigwam} is hexagon 6.
\begin{diagram}
j_{SY!}i_{SQ!}i_{SQ}^*j_{SY}^*\cf & \rTTo^\sim & k_{QY!}k_{QY}^*\cf &
\rTTo^\sim & i_{YZ!}j_{QZ!}j_{QZ}^*i_{YZ}^*\cf \\
\uTTo<{j_{SY!}b_{j^*\cf}^{SQ}} &&&&
\dTTo>{i_{YZ!}a_{i^*\cf}^{QZ}} \\
j_{SY!}j_{SY}^*\cf & \rTTo^{a_\cf^{SY}} & \cf &
\rTTo^{b_\cf^{YX}} & i_{YZ!}i_{YZ}^*\cf
\end{diagram}
Due to identifications
\begin{gather*}
j_{SY!}i_{SQ!}i_{SQ}^*j_{SY}^*\cf = \bigl((\cf\big|_Q)^S\bigr)^Y =
(\cf\big|_Q)^Y = \cf_Q = k_{QY!}k_{QY}^*\cf, \\
i_{YZ!}j_{QZ!}j_{QZ}^*i_{YZ}^*\cf = \bigl((\cf\big|_Q)^Z\bigr)^Y =
(\cf\big|_Q)^Y = \cf_Q = k_{QY!}k_{QY}^*\cf
\end{gather*}
this hexagon is in fact the square
\begin{diagram}
\cf_S & \rMono & \cf \\
\dEpi && \dEpi \\
\cf_Q & \rMono & \cf_Z
\end{diagram}
Considering mappings of stalks of these sheaves at a point of $Y$ we
deduce commutativity of this square.
\end{proof}

\begin{corollary}\label{cor-octa-equiv}
The octahedron of \propref{pro-stan-octa} holds true in
$G$\n-equivariant case.
\end{corollary}

\begin{proof}
The octahedron for closed $G$\n-subspaces $X\subset Z\subset Y$ in
$D_G(Y)$ follows from that for closed subspaces
\( \overline{X\times_YP} \subset \overline{Z\times_YP} \subset
\overline{P} \), where $P$ is a $G$\n-resolution of $Y$,
see \secref{equi-dis-tria}.
\end{proof}

%\tableofcontents
%\bibliographystyle{klunum}
%\bibliography{yuri}
%\end{document}

\end{document}